\documentclass[preprint,12pt]{elsarticle}

\usepackage{amssymb}
\usepackage{amsthm}
\usepackage{subcaption}

\newtheorem{theorem}{Theorem}[section]

\theoremstyle{definition}

\theoremstyle{remark}
\newtheorem{remark}[theorem]{Remark}

\def \bc{\begin{center}}
\def \ec{\end{center}}

\begin{document}

\begin{frontmatter}




\title{Subspace method based on neural networks  for solving the partial differential equation }

\author[mymainaddress,mythirdaddress]{Zhaodong Xu}\ead{ xuzhaodong\_math@163.com}
\author[mymainaddress,mysecondaryaddress]{Zhiqiang Sheng  \corref{cor1}}\ead{sheng\_zhiqiang@iapcm.ac.cn}

\cortext[cor1]{Corresponding author.}

\address[mymainaddress]{National Key Laboratory of Computational Physics, Institute of Applied Physics and Computational Mathematics,  Beijing, 100088, China}
\address[mysecondaryaddress]{HEDPS, Center for Applied Physics and Technology, and College of Engineering, Peking University, Beijing, 100871, China}
\address[mythirdaddress]{Graduate School of China Academy of Engineering Physics, Beijing 100088, China}

\begin{abstract}
We  present a subspace method based on neural networks (SNN) for solving the partial differential equation with high accuracy. 
The basic idea of our method is to use some functions based on neural networks  as base functions to span a subspace, then find an approximate solution in this subspace.
We design two special algorithms in  the strong form of partial differential equation.
One algorithm enforces the equation and initial boundary conditions to hold on some collocation points,
and another algorithm enforces $L^2$-norm of the residual of the equation and  initial boundary conditions to be $0$.
Our method can achieve high accuracy with low cost of training. Moreover, our method is free of parameters that need to be artificially adjusted.
Numerical examples show that the cost of training these base functions of subspace
is low, and only one hundred to two thousand epochs are needed for most  tests.
The error of our method can even fall below the level of $10^{-10}$ for some tests.
The performance of our method significantly surpasses the performance of PINN and DGM in terms of the accuracy and computational cost.

\end{abstract}


\end{frontmatter}

\date{}



\section{Introduction}\label{introduce}

Due to the rapid development of machine learning methods, 
the  method based on neural networks attracts more and more attention.
Since neural networks can be used to approximate any function,
they can be used to  approximate the  solution of partial differential equation(PDE).
Researchers have proposed many  numerical methods based on neural networks,
which can be a promising approach  for solving  the partial differential equation.

Many machine learning methods for solving the partial differential equation are based on deep neural networks,
such as physical information neural networks (PINN)\cite{Raissi-19-jcp}, deep Galerkin method (DGM)\cite{Sirignano-18-jcp},
deep Ritz method (DRM)\cite{E-18}, and weak adversarial networks (WAN)\cite{Zang-20-jcp}.
The main difference between these methods is in the construction of  the loss function.
For example, the loss functions of PINN  and DGM are based on the $L^2$-norm
of the residuals of partial differential equation in  strong form.
The loss function of the deep Ritz method is  based on the energy functional corresponding to the weak form of partial differential equation.
The weak adversarial networks method constructs a loss function by minimizing an operator norm
induced from the weak form of partial differential equation.

A  method using deep learning approach for interface problem is proposed in \cite{Wang-20-jcp}.
They   reformulate the equation in variational form, use deep neural networks to represent the solution of the equation,
and use a shallow neural networks to approximate inhomogeneous boundary conditions.
Another approach to solve this problem is proposed in \cite{He-22-jcam},
utilizing different neural networks  in different sub-domains since the solution may change
dramatically across the interface.
A multi-scale fusion networks is constructed in \cite{Ying-23-amm}.
This multi-scale fusion networks can better capture discontinuity, thereby improving accuracy.
A deep learning method based on PINN for multi-medium diffusion problem is proposed in \cite{Yao-23-amame}.
They add the interface continuity condition as a loss term to the loss function
and propose a domain separation strategy.
A cusp-capturing PINN  to solve  interface problem is proposed in \cite{Tseng-23-jcp},
which introduces a cusp-enforced level set function to the networks.

The PINN method is used to solve  one dimensional and two dimensional Euler equations that model high-speed aerodynamic flows in \cite{Mao-20-cmame}.
A conservative PINN on discrete sub-domains for nonlinear conservation laws is proposed in \cite{Jagtap-20-cmame}.
Thermodynamically consistent PINN for hyperbolic systems is presented in  \cite{Patel-22-jcp}.
PINN is used to solve the inverse problems in supersonic flows in \cite{Jagtap-22-jcp}.
A PINN method with equation weight   is introduced in \cite{Liu-24-jsc},
which introduces a weight such that the neural networks concentrate on training the smoother parts of the solutions.

Due to the curse of dimensionality,
deep neural networks methods have been widely applied for solving high-dimensional partial differential equations.
A type of tensor neural networks is introduced in \cite{Wang-doctor,Wang02754}.
They develop an efficient numerical integration method for the functions of the tensor neural network,
and prove the computational complexity to be the polynomial scale of the dimension.
A machine learning method  solving high-dimensional partial differential equation
by using tensor neural networks and a posteriori error estimator is proposed in \cite{Wang02732}.
They use a posteriori error estimator as the loss function to update
these parameters of tensor neural networks.

In recent years, many methods based on shallow neural networks have also received attention,
such as  methods based on extreme learning machine (ELM)\cite{Huang-06} and random feature methods.
ELM-based methods for solving the partial differential equation
have been developed \cite{Dong-21-cmame,Dong-21-jcp, Dong-22-amame,Dong-22-jcp,Li-03087,Ni-23-jsc,Shang-23,Sun-24-jcam}.
A numerical method for solving linear and nonlinear partial differential equations based on neural networks by combining the ideas of extreme learning machines, domain decomposition, and local neural networks  is proposed in \cite{Dong-21-cmame}.
The weight/bias coefficients of all hidden layers in the local neural networks are all preset random values in the interval $[-R_m,R_m]$,
where $R_m$ is a hyperparameter,
and only the weight coefficients of the output layer need to be solved by the least squares method.
A modified batch intrinsic plasticity  method for pre-training the random coefficients  is proposed in \cite{Dong-21-jcp}
in order to reduce the impact of the hyperparameter on accuracy.
A method based on the differential evolution algorithm to calculate the optimal or
near-optimal value of  the hyperparameter is given in \cite{Dong-22-jcp}.
An approach for solving the partial differential equation  based on randomized neural networks and the Petrov-Galerkin method is proposed in \cite{Shang-23}.
 They  allow for a flexible choice of test functions, such as finite element basis functions.
A local randomized neural networks method with discontinuous Galerkin  methods for partial differential equation is developed in \cite{Sun-24-jcam},
which  uses randomized neural networks  to approximate the solutions on sub-domains,
and uses the discontinuous Galerkin  method to glue them together.
A local randomized neural networks method for interface problems is developed in \cite{Li-03087}.
A discontinuous capturing shallow neural networks method for the elliptical interface problem is developed in \cite{Hu-22-jcp}.

The random feature method for solving the partial differential equation is proposed in \cite{Chen-22-jml}.
This method is a natural bridge between traditional and machine learning-based algorithms.
They use  random feature functions to  approximate the solution,  collocation method to take care of the partial differential equation,
and  penalty method to treat the boundary conditions.
A neural networks method which automatically satisfies boundary and initial conditions is proposed in \cite{Lyu-21-csiam}.
A deep mixed residual method  for solving the partial differential equation with high-order derivatives is proposed in \cite{Lyu-22-jcp}.
They  rewrite a high-order partial differential equation into a first-order system, and
use the residual of first-order system  as the loss function.
A random feature method for solving interface problem is proposed in \cite{Chi-24-cmame},
which utilizes two sets of random feature functions on each side of the interface.

Although deep neural networks-based methods have achieved significant progress in solving the partial differential equation,
they suffer from some limitations.
The first limitation is that the accuracy of  these methods is unsatisfactory.
A survey of related literatures shows that the  error of most deep neural networks-based  methods is difficult to fall below the level of $10^{-4}$. Increasing  number of training epochs does not significantly reduce the error.
Another limitation is low  efficiency.
The computational cost of solving the partial differential equation with these methods based on deep neural networks is extremely high.
A lot of computational time is needed for training.
For example, some methods  based on deep neural networks need several hours to achieve certain accuracy,
while traditional methods such as finite element methods can achieve  similar accuracy in just a few seconds.
Due to low accuracy and high computational cost,
it is a challenge for these methods based on deep neural networks  to compete with traditional methods for low dimensional problems.

The hyperparameter of ELM-based methods has a significant impact on accuracy.
The method with the optimal hyperparameter can achieve high accuracy,
however, the method with an inappropriate hyperparameter results in very poor accuracy.
Selecting an optimal hyperparameter is a challenging problem.

In this paper, we  present a subspace method based on neural networks  for solving the partial differential equation with high accuracy.
The basic idea of our method is to use some functions based on neural networks  as base functions to span a subspace, then find an approximate solution in this subspace.
Our method includes three steps.
First, we  give the neural networks architecture which includes input layer, hidden layer, subspace layer and output layer.
Second, we train these base functions of subspace such that the subspace has effective approximate capability to the solution space of equation.
Third, we find an approximate solution in the subspace to approximate  the solution of the equation.
We design two special algorithms in  the strong form of partial differential equation.
One algorithm enforces the equation and initial boundary conditions to hold on some collocation points, we call this algorithm as SNN in discrete form(SNN-D).
Another algorithm enforces $L^2$-norm of the residual of the equation and   initial boundary conditions
to be $0$, we call this algorithm as SNN in integral form(SNN-I).
Our method can achieve high accuracy, and the cost of training is low.

Our method is free of parameters (including hyperparameter and penalty parameter) that need to be artificially adjusted.
Different from ELM, we do not introduce the hyperparameter since we train these parameters of neural networks.
Different from PINN and DGM, we do not use the initial boundary conditions in the loss function, hence we do not need to introduce the penalty parameter.
When the number of hidden layer reduces to $0$ and the number of training epochs becomes $0$, our method degenerates into ELM.
When we use the loss function including both the PDE loss term and initial boundary loss term in our method, and omit the third step of our method, that is the least squares method is not be used to update these parameters, with the training epochs matching those of PINN or DGM, our method degenerates into PINN or DGM.
We use the Adam method to update neural network parameters, while other methods are viable too.
Additionally, we need to solve an algebraic system  by the  least squares method.

Numerical examples show that the cost of training these base functions of subspace
is low, and only one hundred to two thousand epochs are needed for most tests.
The error of our method can even fall below the level of $10^{-10}$ for some tests.
In general, the accuracy of SNN-D is  higher than that of  SNN-I.
Furthermore, the performance of our method significantly surpasses the performance of PINN and DGM in terms of the accuracy and computational cost.

The remainder of this paper is organized as follows.
In section 2, we   describe  the subspace method based on neural networks for solving the partial differential equation.
In section 3, we present some numerical examples  to test the performance of our method.   At last, we give some conclusions.

\section{Subspace method based on neural networks} 
Consider the  following equation:
 \begin{eqnarray}
\mathcal{A}u(\bold{x})&=&f(\bold{x})  \qquad \ \ \mbox{in}\ \ \Omega,
\label{deq1} \\
\mathcal{B} u(\bold{x}) &=&g(\bold{x})  \qquad \ \  \mbox{on} \ \
\partial \Omega,  \label{deq2}
\end{eqnarray}
where $\bold{x}=(x_1,x_2,\cdots, x_d)^T$,   $\Omega$ is a bounded domain in ${\textrm{R}}^d$, $\partial \Omega$ is the boundary of $\Omega$, $\mathcal{A}$ and $\mathcal{B}$ are the differential operators, $f$ and $g$ are given functions.

\subsection{Neural networks architecture}
In this section, we describe the neural networks architecture. For simplicity,
we only describe the neural networks architecture with  one-dimensional output.
Of course, this neural networks architecture can be used in the case with $k$-dimensional output.

\begin{figure}[htbp]
\begin{center}\includegraphics[width = 10cm]{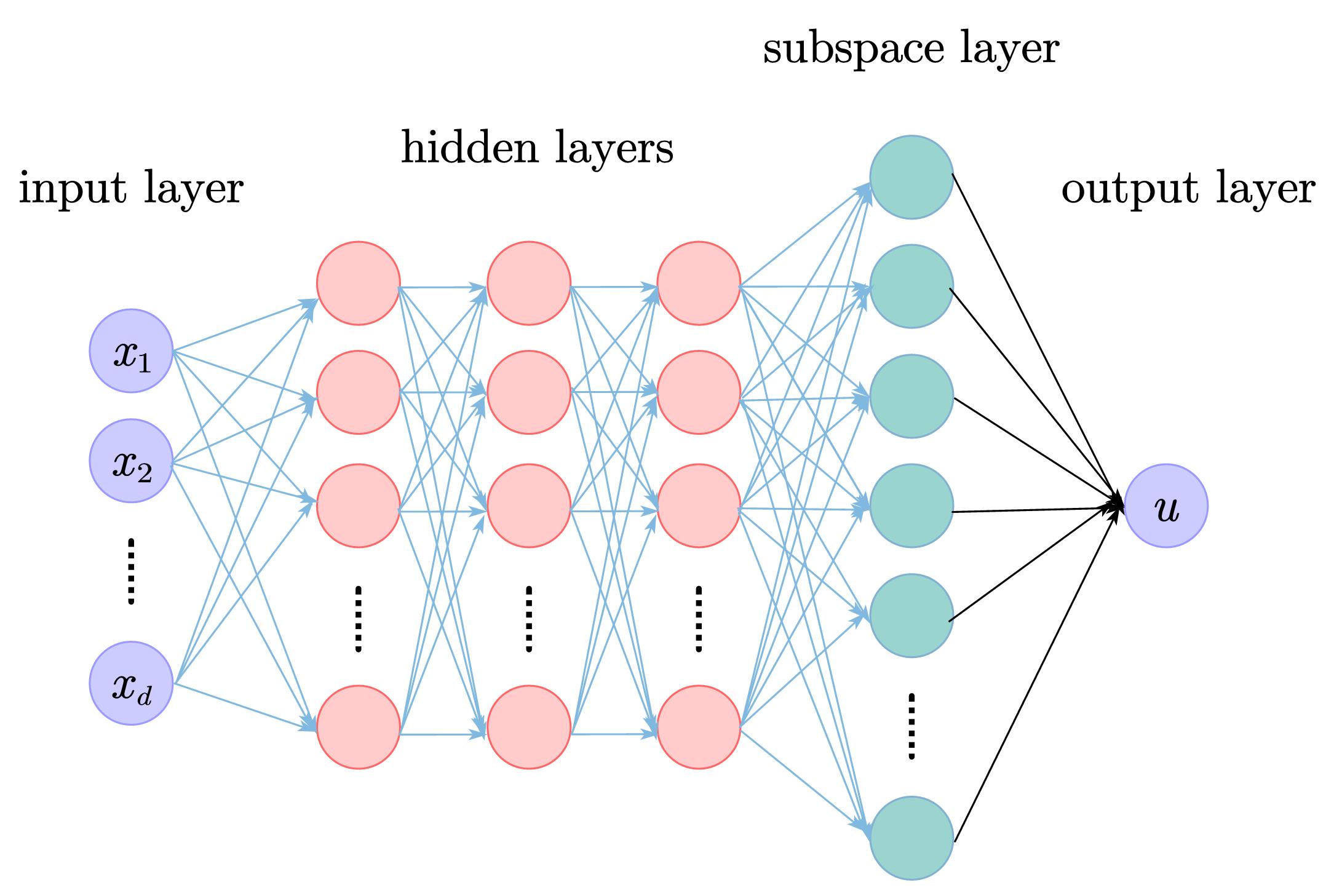}
 \caption{The neural networks architecture.} \label{network}
\end{center}
\end{figure}

The neural networks architecture consists of four key components, including an input layer, hidden layers, a subspace layer and  an output layer. The neural network employs some hidden layers to enrich the expressive capability of network. The subspace layer is essential for constructing a finite dimension space that approximates the solution space of the equation.  Figure \ref{network} illustrates this specialized architecture.

Let $K$ be the number of  hidden layers, $n_1$, $n_2$, $\cdots$, $n_K$ be the number of  neurons in each hidden layer, respectively.
Let $M$ be the dimension of subspace in the subspace layer, and $\varphi_j$($j=1,2,\cdots,M$) be base functions of subspace,
 and $\omega_j$($j=1,2,\cdots,M$) be some coefficients related to base functions.
 Denote $\varphi = (\varphi_1,\varphi_2,\cdots, \varphi_M)^T$ and $\omega =(\omega_1, \omega_2, \cdots, \omega_M)^T$.

The propagation process  can be expressed as follows:
\[
\left\{
\begin{array}{l}
\bold{y}_0=\bold{x}, \\
\bold{y}_k=\sigma_k(W_k \cdot \bold{y}_{k-1} + \bold{b}_k), \ \hfill k =1,2,\cdots, K+1,\\
 \bold{\varphi}=\bold{y}_{K+1},\\
u= \bold{\varphi}\cdot \bold{\omega},
\end{array}\right. 
\]
where $W_k \in {\textrm{R}}^{n_k \times n_{k-1}}$  and $b_k\in {\textrm{R}}^{n_k}$ are the weight and bias, respectively, $n_0=d$ is the dimension of input and
$n_{K+1}=M $ is the dimension of subspace.  $\bold{x} \in {\textrm{R}}^d$ is the input and $\sigma(\cdot)$ is  the activation function.
$\theta =\{  W_1,\cdots, W_{K+1},b_1,\cdots, b_{K+1}\}$ is the set of  parameters in neural networks, $u(x;\theta, \omega)$ is the output with respect to input $x$ with parameters $\theta$ and $\omega$.

These weight and bias coefficients are initially randomly generated and subsequently updated by minimizing the loss function   $\mathcal{L}(\bold{x};\theta,\omega)$.   Usually, this update is achieved by gradient descent method.
In each iteration, these parameters can be updated as follows:
\begin{eqnarray*}
W_k &\leftarrow & W_k -\eta \frac{\partial \mathcal{L}(\bold{x};\theta,\omega)}{\partial W_k},\\
b_k  &\leftarrow & b_k -\eta \frac{\partial \mathcal{L}(\bold{x};\theta,\omega)}{\partial b_k},
\end{eqnarray*}
where $\eta > 0$ is the learning rate.
For the gradient descent method, it is needed to calculate the partial derivative of the loss function with respect to network parameters,
these partial derivatives are implemented through an automatic differentiation mechanism in  Pytorch and Tensorflow.

After the training process, $\bold{\omega}$ is determined by enforcing the equation and boundary conditions to hold.
Although $\bold{\omega}$ can be updated during training,
this is not necessary as its final value is obtained by satisfying these constraints. In fact, we find that there has no significant impact on the numerical results.

\subsection{A general frame of subspace method based on neural networks}
Now, we describe a general frame for the subspace method based on neural networks for solving the partial differential equation.

First, we  construct the neural network architecture which includes input layer, hidden layer, subspace layer and output layer.
Then, we train the base functions of subspace such that the subspace has effective approximate capability to the solution space of equation.
At last, we find an approximate solution in the subspace to approximate  the solution of the equation.
This general frame of SNN is as follows:
 
\vspace{3mm}
\begin{tabular}{l}\hline
 A general frame of SNN\\
\hline
1. Initialize nerual networks architecture, generate randomly $\theta$, and\\
  \ \quad  give $\omega$.\\

2. Update parameter $\theta$ by minimizing the loss function, i.e. training\\
\ \quad    the base functions of the subspace layer $\varphi_1$, $\varphi_2$, $\cdots$, $\varphi_M$. \\ 

3. Update parameter $\omega$ and find an approximate solution in the\\
  \ \quad  subspace to approximate  the  solution of the equation.  \\
\hline
\end{tabular}
\vspace{3mm}

\begin{remark}
Many studies have shown that the imbalance between PDE loss and initial boundary loss in the training process can lead to
lose the accuracy and increase significantly the cost of training.
Weighted techniques have been used to correct this imbalance\cite{McClenny-23-jcp, Wight-21-cicp, Yao-23}.
However, how to determine these weights is a challenging problem.
We find that the initial boundary loss is not important in the training process of base functions for many problems.
In order to overcome this challenge  of selecting weights (penalty parameter), we do not use the initial boundary conditions in step 2, hence we do not
need to introduce the penalty parameter. In fact, the information of PDE is enough to train the base functions of subspace  for many problems,
and the information of  initial boundary conditions is not necessary.
Of course, if one does not care about the cost of training, the loss function including both PDE loss and initial boundary loss can also be used in our method.

\end{remark}

\begin{remark}
In step 2, we fix the parameter $\omega$, and train the parameter $\theta$ by  minimizing the loss function.
The  aim is to derive suitable base functions such that the subspace spanned by these base functions 
has effective approximate properties. In order to make a balance for the accuracy and efficiency, it is not necessary to solve
 the minimization problem  accurately.
\end{remark}

\begin{remark}
In step 3, in order to find an approximate solution of the equation, we need to use the information of both PDE and initial boundary conditions.
We can obtain an algebraic system, and solve this system to get $\omega$.
In general, this algebraic system does not form a square matrix, it is typically solved using the least squares method.
\end{remark}

\subsection{Two special algorithms}
In this section, we will describe two special algorithms.
We only consider the strong form of partial differential equation in this paper.
One algorithm enforces the equation and initial boundary conditions to hold on some collocation points. We call
this algorithm as the subspace method based on neural networks in discrete form (SNN-D).
Another algorithm enforces the $L^2$-norm of the residual of the equation and the  $L^2$-norm of the residual of the initial boundary conditions
to be $0$. We call this algorithm as the subspace method based on neural networks in integral form (SNN-I).

\subsection{The algorithm in discrete form: SNN-D}
Denote the coordinates of inner collocation points in $\Omega$ by $\bold{x}_j$,
the coordinates of boundary collocation points on $\partial \Omega$ by $\bold{\bar{x}}_j$.
In general, the loss function is as follows:
\begin{equation}
\mathcal{L}=\frac1N \sum_{j=1}^N \left(\mathcal{A} u(\bold{x}_j;\theta,\omega) -f(\bold{x}_j)\right)^2
+\lambda \frac{1} {\bar{N}} \sum_{j=1}^{\bar{N}} \left(\mathcal{B} u(\bold{\bar{x}}_j;\theta,\omega) -g(\bold{\bar{x}}_j)\right)^2,
\label{loss-g}
\end{equation}
where $N$ is the number   of inner collocation points in $\Omega$, $\bar{N}$ is the number   of boundary collocation points on $\partial \Omega$,
and $\lambda$ is a penalty parameter. In order to make a balance for the loss between different loss terms in the training process,
one needs to select an optimal parameter $\lambda$.
However, it is a challenging problem to select an optimal parameter $\lambda$.
Many literatures have focused on how to give an optimal parameter.

To avoid introducing this penalty parameter,  we  define the following loss function which  contains only the PDE loss term. 
\begin{equation}
\mathcal{L_D}=\frac1N \sum_{j=1}^N \left(\mathcal{A} u(\bold{x}_j;\theta,\omega) -f(\bold{x}_j)\right)^2.
\label{loss-d}
\end{equation}
It is obvious  that this loss function  includes only the information of PDE itself, and omit the information of initial boundary conditions.

First, we fix the parameter $\omega$, and train the parameter $\theta$ by  minimizing the loss function $\mathcal{L_D}(\bold{x};\theta,\omega)$ in (\ref{loss-d}).
For example, we can take $\omega_j=1$($j=1,\cdots, M$), which can ensure that each base function has the same contribution to the loss function
in the training process.
In this step, our main aim is to obtain some suitable base functions of subspace, and expect that this subspace
has  effective approximate properties.
In fact, we find that there has no significant impact on the the numerical results for choosing  different $\omega$.
We use Adam method to solve this minimization problem.
In order to balance the accuracy and efficiency, it is not necessary to solve
 the minimization problem accurately. Let $\mathcal{L_{D}}_0(\bold{x};\theta,\omega)$ be the initial loss.
 The training process stops if the following condition is satisfied,
\begin{equation}
\frac{\mathcal{L_D}(\bold{x};\theta,\omega)}{\mathcal{L_{D}}_0(\bold{x};\theta,\omega)} \le \varepsilon.
\label{loss-dstop}
\end{equation}
This implies that  the training process stops if the  loss decreases to a certain level.
To prevent excessive  training epochs, we introduce the maximum number of training epochs $N_{max}$.
If the number of training epochs reaches  $N_{max}$,  the training process also stops.
Since it is not necessary to solve
 the minimization problem  accurately, $\varepsilon$ does not need to be too small and $N_{max}$ does not need to be too big.
 For example, we  take $\varepsilon=10^{-3}$ and $N_{max}=5000$. For most tests,  only   one hundred to two thousand epochs are needed.

Then, we fix the parameter $\theta$ and update the parameter $\omega$. Let
\begin{equation}
u(\bold{x}) = \sum_{j=1}^M \omega_j \varphi_j(\bold{x}).
\label{u}
\end{equation}
Substitute the expression above into the equation (\ref{deq1}) and (\ref{deq2}) to obtain
 \begin{eqnarray}
\sum_{j=1}^M \omega_j \mathcal{A}\varphi_j(\bold{x})&=&f(\bold{x}),  \label{Au} \\
\sum_{j=1}^M \omega_j \mathcal{B}\varphi_j(\bold{x}) &=&g(\bold{x}).  \label{Bu}
\end{eqnarray}

We enforce equation (\ref{Au}) on all inner collocation points $\bold{x}_j$($j=1,2,\cdots, N$), and enforce equation (\ref{Bu}) on all boundary
collocation points $\bold{\bar{x}}_j$($j=1,2,\cdots, \bar{N}$), which lead to the following equations,
 \begin{eqnarray}
\sum_{j=1}^M \omega_j \mathcal{A}\varphi_j(\bold{x}_i)&=&f(\bold{x}_i),  \quad i=1,2,\cdots,N, \label{Au-d} \\
\sum_{j=1}^M \omega_j \mathcal{B}\varphi_j(\bold{\bar{x}}_i) &=&g(\bold{\bar{x}}_i), \quad i=1,2,\cdots, \bar{N}.  \label{Bu-d}
\end{eqnarray}
Hence, we can obtain an algebraic  system, which consists of $N+\bar{N}$ equations, and  $M$ unknowns $\omega_1, \omega_2,\cdots, \omega_M$.
In this algebraic system, $\mathcal{A}\varphi_j(\bold{x}_i)$ and $\mathcal{B}\varphi_j(\bold{\bar{x}}_i)$ are all known.
We use the least squares method to solve this system.
At last, we can obtain an approximate solution $u$.

We summarize the main steps of SNN-D in Algorithm 1.

\begin{tabular}{l}\hline
 {\bf Algorithm 1:}  SNN-D\\
\hline
1. Initialize nerual networks architecture.\\

2. Generate randomly $\theta$, and give $\omega$.\\

3. Update parameter $\theta$ by minimizing the loss function $\mathcal{L_D}(\bold{x};\theta,\omega)$  \\
   \ \quad in (\ref{loss-d}) until (\ref{loss-dstop}) holds or the epochs reach $N_{max}$. \\

4. Obtain the base functions of the subspace $\varphi_1$, $\varphi_2$, $\cdots$, $\varphi_M$.\\

5. Solve the algebraic system resulted from (\ref{Au-d}) and (\ref{Bu-d}) to update $\omega$.\\

6. Obtain an approximate solution $u$.  \\
\hline

\end{tabular}

\vspace{0.3cm}

\begin{remark}
When the number of hidden layer reduces to $0$ and the number of training epochs becomes $0$, SNN-D degenerates into ELM.
When we use the loss function including both the PDE loss term and initial boundary loss term, and omit steps 4 and 5 in  Algorithm 1, with the training epochs matching those  of PINN, SNN-D degenerates into PINN.
\end{remark}

\subsection{The algorithm in integral form: SNN-I}
Define  the loss function as follows:
\begin{equation}
\mathcal{L_I}=\| \mathcal{A} u(\bold{x};\theta,\omega) -f(\bold{x})\|_{L^2(\Omega)}^2.
\label{loss-i}
\end{equation}
We can see that the loss function includes  only the information of PDE, and omit the information of initial boundary conditions.

At first, we fix the parameter $\omega$, and train the parameter $\theta$ by  minimizing the loss function $\mathcal{L_I}(\bold{x};\theta,\omega)$ in (\ref{loss-i}).
Similar to SNN-D, we can take $\omega_j=1$($j=1,\cdots, M$), which can ensure that each base function has the same contribution to the loss function
in the training process.  Let $\mathcal{L_{I}}_0(\bold{x};\theta,\omega)$ be the initial loss.
 The training process stops if the following condition is satisfied,
\begin{equation}
\frac{\mathcal{L_I}(\bold{x};\theta,\omega)}{\mathcal{L_{I}}_0(\bold{x};\theta,\omega)} \le \varepsilon.
\label{loss-istop}
\end{equation}
This implies that  the training process stops if the  loss decreases to a certain level.
If the epochs reach $N_{max}$,  the training process also stops.
Similar to SNN-D,  it is not necessary to solve  the minimization problem accurately in order to balance the accuracy and efficiency.
In this step, our main aim is to obtain some suitable base functions of subspace, and expect that this subspace has  effective approximate properties.

Then, we fix the parameter $\theta$ and update the parameter $\omega$.
We can obtain the expression of $u$ (see (\ref{u})) after we obtain base functions of subspace.
Let
\begin{eqnarray}
F_1 (\omega_1, \omega_2, \cdots, \omega_M)&=& \| \mathcal{A} u(\bold{x};\theta,\omega) -f(\bold{x})\|_{L^2(\Omega)}^2,\label{F1-0}\\
F_2 (\omega_1, \omega_2, \cdots, \omega_M)&=& \| \mathcal{B} u(\bold{x};\theta,\omega) -g(\bold{x})\|_{L^2(\partial \Omega)}^2. \label{F2-0}
\end{eqnarray}
Substitute the expression (\ref{u}) into (\ref{F1-0}) and (\ref{F2-0}) to obtain
\begin{eqnarray}
F_1& = &\| \sum_{j=1}^M \omega_j  \mathcal{A}\varphi_j -f\|_{L^2(\Omega)}^2
=\int_{\Omega} \left(\sum_{j=1}^M \omega_j \mathcal{A}\varphi_j-f \right)^2 d\bold{x},\label{F1}\\
F_2 &=& \|  \sum_{j=1}^M \omega_j \mathcal{B} \varphi_j -g\|_{L^2(\partial \Omega)}^2
=\int_{\partial \Omega} \left(\sum_{j=1}^M \omega_j \mathcal{B}\varphi_j-g \right)^2 ds. \label{F2}
\end{eqnarray}
From the definitions of $F_1$ and $F_2$, we can see there are $F_1 \ge 0$ and $F_2 \ge 0$.
From (\ref{F1}) and (\ref{F2}), we can know that $F_1$ and $F_2$ are a quadratic function about $\omega_1, \omega_2, \cdots, \omega_M$.
In order to minimize $F_1$ and $F_2$, we let
\begin{eqnarray}
\frac{\partial F_1(\omega_1, \omega_2, \cdots, \omega_M)}{\partial \omega_i}  &=& 0, \quad i=1,2,\cdots, M,\\
\frac{\partial F_2(\omega_1, \omega_2, \cdots, \omega_M)}{\partial \omega_i}  &=& 0, \quad i=1,2,\cdots, M.
\end{eqnarray}
Which leads to
\begin{eqnarray}
\sum_{j=1}^M \left( \int_{\Omega} \mathcal{A}\varphi_i \cdot \mathcal{A}\varphi_j d\bold{x}\right) \omega_j&=& \int_{\Omega} \mathcal{A}\varphi_i \cdot f d\bold{x}, \quad i=1,2,\cdots, M,\label{Au-i}\\
\sum_{j=1}^M \left( \int_{\partial \Omega} \mathcal{B}\varphi_i \cdot \mathcal{B}\varphi_j ds\right) \omega_j&=& \int_{\partial \Omega} \mathcal{B}\varphi_i \cdot g ds, \quad i=1,2,\cdots, M.\label{Bu-i}
\end{eqnarray}
Hence, we can obtain an algebraic system about $\omega_1, \omega_2, \cdots, \omega_M$. There has 2$M$ equations and $M$ unknowns,
we use the least squares method to solve it.
At last, we can obtain an approximate solution $u$.

We summarize the main steps of SNN-I in Algorithm 2.

\vspace{3mm}
\begin{tabular}{l}\hline
 {\bf Algorithm 2:}  SNN-I\\
\hline
1. Initialize neural networks architecture.\\

2. Generate randomly $\theta$, and give $\omega$.\\

3. Update parameter $\theta$ by minimizing the loss function $\mathcal{L_I}(\bold{x};\theta,\omega)$  \\
  \ \quad  in (\ref{loss-i}) until (\ref{loss-istop}) holds or the epochs reach $N_{max}$. \\

4. Obtain the base functions of the subspace $\varphi_1$, $\varphi_2$, $\cdots$, $\varphi_M$.\\

5. Solve the algebraic system resulted from (\ref{Au-i}) and (\ref{Bu-i}) to update $\omega$.\\

6. Obtain an approximate solution $u$.  \\
\hline
\end{tabular}
\vspace{3mm}

\begin{remark}
When we use the loss function including both the PDE loss term and initial boundary loss term, and omit steps 4 and 5 in  Algorithm 2, with the training epochs matching those  of DGM, SNN-I degenerates into DGM.
\end{remark}

\section{Numerical results}
\label{Numerical results}
In this section, we present several numerical experiments to demonstrate the performance of our method. In Section \ref {Helmholtz equation}, we test the performance of SNN by solving the one-dimensional Helmholtz equation. We give some numerical results with various network depth and subspace dimension. In Section \ref {Poisson equation}, we show the results of SNN for solving the two-dimensional Poisson equation. In Section \ref {Advection equation}, we test SNN for solving the advection equation. Specifically, we demonstrate the impact of initial boundary conditions on  the training of base functions. In Section \ref {Parabolic equation}, we test the performance of SNN for solving the parabolic equation. In Section \ref{Anisotropic diffusion equation}, we use SNN to solve the strongly anisotropic diffusion equation.

We employ the deep learning framework PyTorch \cite{paszke-19-ansps} for code development, ensuring all variable data types are set to float64. The numerical results may vary due to differences in model architecture, the size of the training data, weights, optimizers, etc. To demonstrate  the robustness and accuracy of our method, we maintain consistent settings for all numerical examples in this paper. Specifically, our numerical experiments utilize a feedforward fully connected neural network (FNN) with four hidden layers, each containing 100 neurons. The activation function is the Tanh function, and the subspace dimension is uniformly fixed at 300.

In all numerical experiments, we use the Adam optimizer. The training process stops  when the relative  loss is less than $\varepsilon$ or  the epochs reach   $N_{max}$. In all tests, we take $\varepsilon=1e-3$ and $N_{max}=5000$. The settings of optimizer , including the learning rate, are kept at their default values. Neural network parameters are randomly generated by using the Xavier method.  For the coefficient matrix, a preprocessing step is employed to ensure the maximum value in each  row is normalized to 1. The right hand side is  processed accordingly.

We compare our method with the existing methods, i.e, PINN, DGM and ELM.
For comparison purposes, identical network architectures, parameter settings, and initialization methods are employed for both PINN and DGM, with the Adam optimizer training through 50000 epochs. PINN utilizes the same sampling points as SNN-D, while DGM employs integration points identical to those in SNN-I. In order to accurately reproduce numerical results, the seed for generating random numbers  in all numerical experiments is set to 1. For ELM, the hyperparameter $R_m=1$ is specified, meaning all network parameters are randomly generated within the range of $[-1,1]$, and the  hidden layer contains 300 neurons, matching the dimension of subspace.

For SNN-D, PINN and ELM, we evaluate the accuracy  by using the point-wise relative $L^2$ error, defined as follows:
\begin{eqnarray*}
\|e\|_{ {L}^2}=\frac{\sqrt{\sum_{i=1}^N\left|u_\theta\left(X_i\right)
-u^*\left(X_i\right)\right|^2}}{\sqrt{\sum_{i=1}^N\left|u^*\left(X_i\right)\right|^2}}.
\end{eqnarray*}
For SNN-I and DGM, we evaluate the accuracy by using the relative $L^2$ error in integral form, defined as follows:
\begin{eqnarray*}
\|e\|_{ {L}^2} = \frac{\sqrt{\int_{\Omega}\left|u_\theta(x) - u^*(x)\right|^2 \,dx}}{\sqrt{\int_{\Omega}\left|u^*(x)\right|^2 \,dx}},
\end{eqnarray*}
where $ u^*$ is the exact solution, and $u_\theta$ is the approximate solution.

\subsection{Helmholtz equation}
\label{Helmholtz equation}
In the first test, we consider the one-dimensional  Helmholtz equation  on
the domain $\Omega = (a, b)$, see \cite{Dong-21-cmame},
\begin{eqnarray}
\label{deq23}
\frac{d^2u}{dx^2} - \lambda u = f(x),
\end{eqnarray}
with the boundary conditions $u(a) = h_1$ and $u(b) = h_2$, where  $f (x)$ is the source term, and $h_1$ and $h_2$ are the boundary conditions. In this test, we take $\lambda = 10$, $a = 0$ and $b = 2$. We
choose the source term $f (x)$ such that Eq. (\ref{deq23}) has the following solution,
\begin{eqnarray*}
u(x) = \sin\left(3\pi x + \frac{3\pi}{20}\right) \cos\left(2\pi x + \frac{\pi}{10}\right) + 2.
\end{eqnarray*}

We use SNN method to solve  this problem. For SNN-D, 1000 points are equally sampled across the interval [a, b]. For SNN-I, the complex Gaussian quadrature formula is applied to segment [a, b] into 30 sub-intervals, with each sub-interval containing 10 points.

\begin{figure}[htbp]
\begin{center}\includegraphics[width = 6cm]{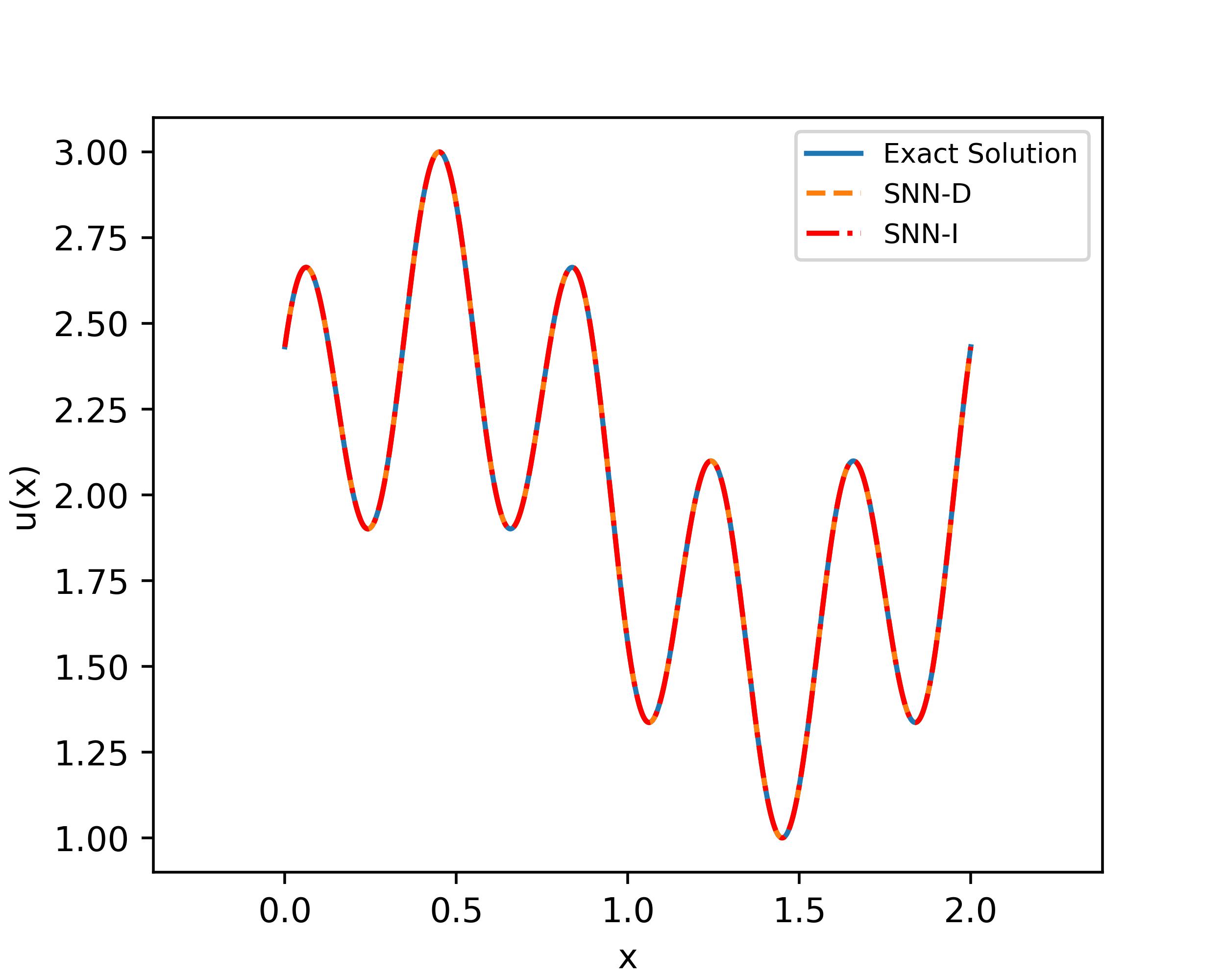}
 \caption{Solution obtained by SNN-D and SNN-I for Helmholtz Equation.} \label{1d_hz_exact_vs_appro}
\end{center}
\end{figure}

\begin{figure}[htbp]
\centering
\subcaptionbox{SNN-D}{\includegraphics[width=0.45\linewidth]{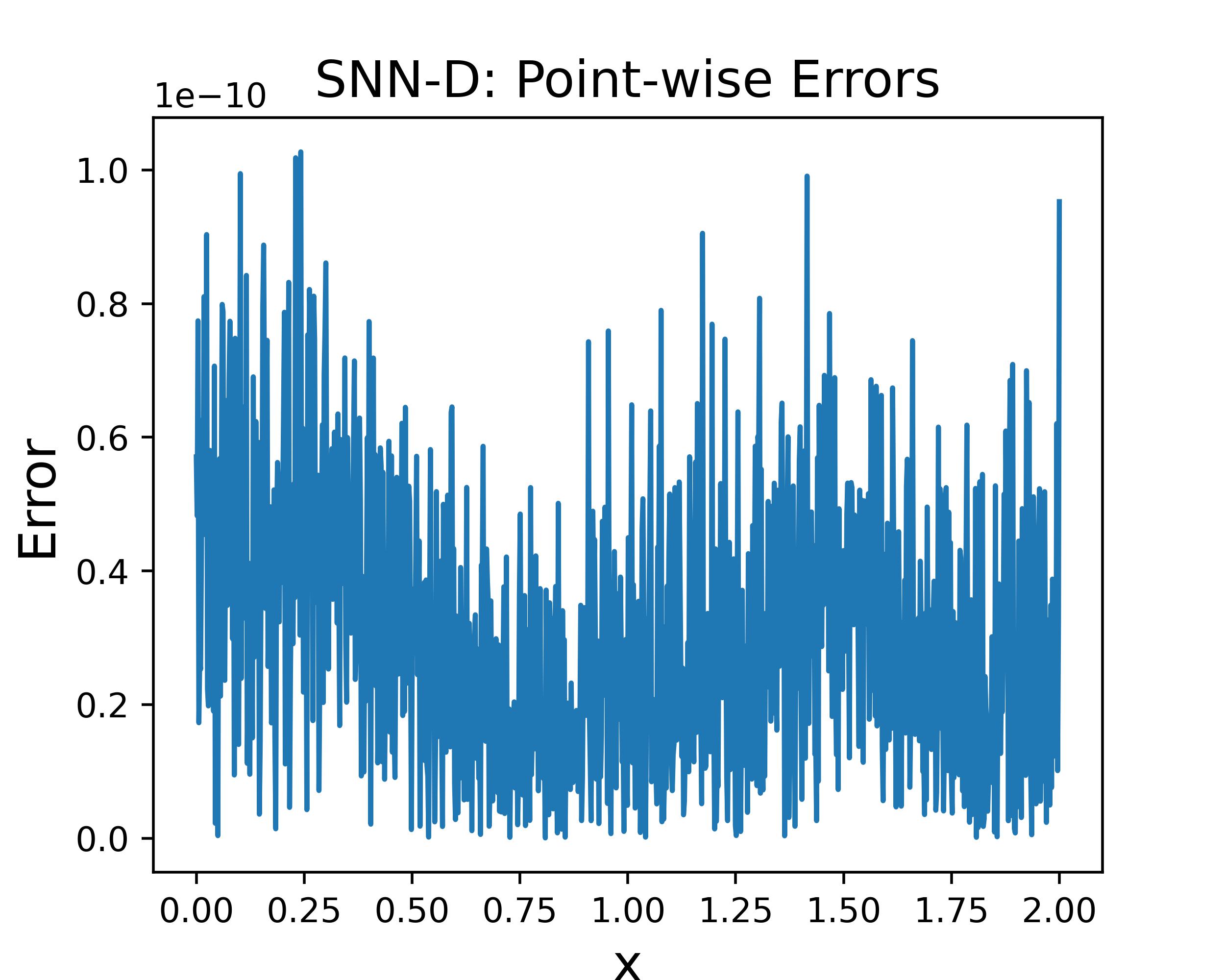}}
\hspace{0.3cm}  
\subcaptionbox{SNN-I}{\includegraphics[width=0.45\linewidth]{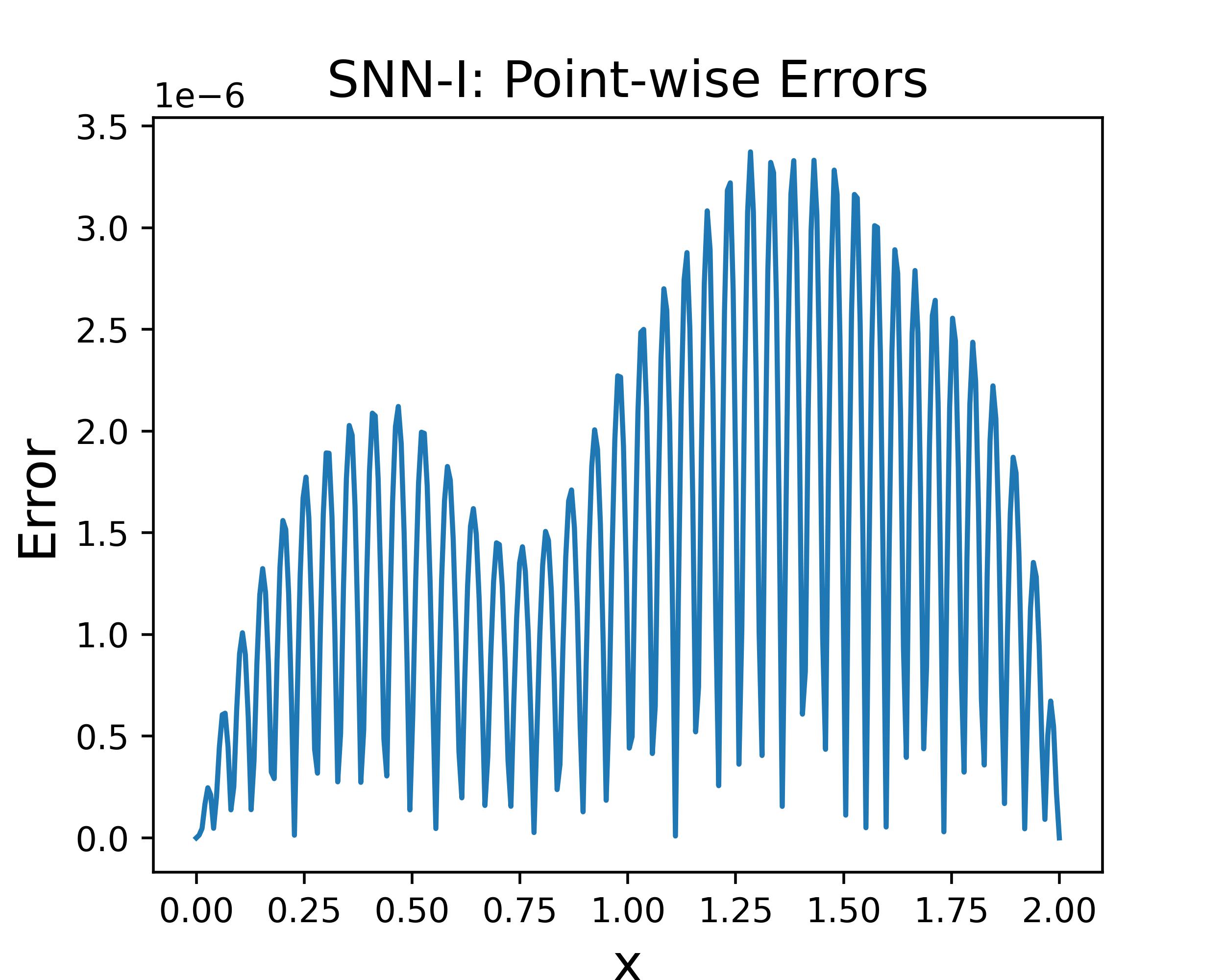}}
\caption{ Point-wise errors of  SNN-D and SNN-I for Helmholtz equation.}
\label{1d_hz_point_error_fig}
\end{figure}

Figure \ref{1d_hz_exact_vs_appro} illustrates the  solutions obtained by our method. 
Figure \ref{1d_hz_point_error_fig} displays the point-wise errors  for both algorithms. Table \ref{1d_hz_error} presents the relative $L^2$ error and $L^\infty$ error for different methods. Notably, after    50000 training epochs, the relative $L^2$ errors for PINN and DGM are 5.69e-03 and 2.29e-04, respectively. 
Increasing the number of training epochs does not significantly reduce the error.
Numerical results demonstrate that  SNN-D  achieves relative $L^2$ error of 1.79e-11 with only 1960  epochs, 
 and  SNN-I achieves relative $L^2$ error of 8.02e-07 with only 3708  epochs. 
This indicates that SNN method can achieves high accuracy with significantly fewer epochs. 
With $R_m=1$, the relative $L^2$ error of ELM is 8.70e-06. Noticely, reducing the number of hidden layer and the number of training epochs to zero will lead SNN-D to degenerate into ELM. This implies that the subspace  consists of some random functions and does not contain the information of PDE.

\begin{table}[!htbp]
\caption{The errors and epochs  of different methods for Helmholtz equation. }\label{1d_hz_error}
\begin{center}
\small
\begin{tabular}{cccccc}
\hline
Method & $\|e\|_{L^2}$ & $\|e\|_{L^\infty}$ & epochs\\
\hline
PINN & 5.69e-03 & 1.99e-02 & 50000 \\
ELM & 8.70e-06 & 6.66e-05 & 0 \\
SNN-D & 1.79e-11 & 1.03e-10 & 1960 \\
DGM & 2.29e-04 & 2.02e-03 & 50000 \\
SNN-I & 8.02e-07 & 3.37e-06 & 3708 \\
\hline
\end{tabular}
\end{center}
\end{table}

Tables \ref{1dhz_SNN_D_point_subspace} and \ref{1dhz_SNN_I_point_subspace} present the numerical results of SNN-D and SNN-I for various numbers of sampling points and subspace dimension. The network setup is uniform with four hidden layers, each with 100 neurons. SNN-D utilizes uniform sampling, whereas SNN-I adopts the composite Gaussian quadrature rule, sampling ten points per subinterval and incrementing subinterval counts for accuracy enhancement. It becomes apparent that the error decreases with increasing both sampling points and subspace dimension. 
With a subspace dimension of 20, the error primarily ranges between $10^ {-2}$ and $10^{-1}$.
Expanding the subspace dimension to 40,  SNN-D can achieve the  accuracy of $10^{-3}$  and SNN-I can achieve the accuracy $10^{-5}$. Upon reaching a subspace dimension of 60, the expressive capacity markedly improves. As the number of sampling points increases,  SNN-D achieves the accuracy of $10^{-8}$ and SNN-I achieves the accuracy of $10^{-7}$. Further increasing the subspace dimension,  the error of SNN-I does not further decrease, while the error of SNN-D decreases to $10^{-11}$. The number of training epochs is from 400 to 2000 for SNN-D, and from 400 to 4000 for SNN-I, significantly less than that required by PINN and DGM. 

\begin{table}[!htbp]
\caption{The errors and epochs of SNN-D across various numbers of sampling points and subspace dimension $M$ for  Helmholtz equation.}\label{1dhz_SNN_D_point_subspace}
\begin{center}
\resizebox{\columnwidth}{!}{%
\small
\begin{tabular}{cccccccc}
\hline
Points & $M$ & 20 & 40 & 60 & 80 & 100 & 300 \\
\hline
20 & $\|e\|_{ {L}^2}$ & 2.13e-01 & 1.49e-02 & 3.19e-03 & 5.63e-03 & 1.25e-02 & 9.88e-03 \\
 & epochs & 710 & 546 & 833 & 760 & 684 & 1170 \\
40 & $\|e\|_{ {L}^2}$ & 2.68e-01 & 7.67e-05 & 4.73e-05 & 2.63e-05 & 2.84e-05 & 1.15e-05 \\
 & epochs & 462 & 457 & 1223 & 813 & 738 & 1514 \\
60 & $\|e\|_{ {L}^2}$ & 3.66e-01 & 3.86e-04 & 3.29e-06 & 1.29e-06 & 2.39e-07 & 7.41e-08 \\
 & epochs & 433 & 659 & 1116 & 685 & 993 & 2065 \\
80 & $\|e\|_{ {L}^2}$ & 3.26e-01 & 1.99e-04 & 1.37e-07 & 1.97e-10 & 1.16e-09 & 3.09e-10 \\
 & epochs & 412 & 802 & 893 & 1450 & 1247 & 1458 \\
100 & $\|e\|_{ {L}^2}$ & 2.68e-01 & 1.60e-03 & 3.53e-08 & 6.72e-09 & 1.01e-10 & 3.03e-11 \\
 & epochs & 421 & 672 & 693 & 1144 & 879 & 1778 \\
300 & $\|e\|_{ {L}^2}$ & 1.68e-01 & 3.78e-03 & 5.72e-08 & 1.13e-10 & 4.20e-11 & 1.38e-10 \\
 & epochs & 465 & 782 & 786 & 1975 & 1614 & 1437 \\
500 & $\|e\|_{ {L}^2}$ & 2.92e-01 & 2.98e-03 & 2.32e-07 & 7.32e-11 & 1.30e-10 & 7.08e-11 \\
 & epochs & 442 & 650 & 699 & 1783 & 1460 & 1966 \\
1000 & $\|e\|_{ {L}^2}$ & 4.73e-01 & 3.83e-02 & 6.79e-08 & 5.91e-11 & 2.58e-11 & 1.79e-11 \\
 & epochs & 404 & 630 & 729 & 1179 & 1184 & 1960 \\
\hline
\end{tabular}
} 
\end{center}
\end{table}

\begin{table}[!htbp]
\caption{The errors and epochs of SNN-I across various numbers of sampling points and subspace dimension M for Helmholtz equation.}\label{1dhz_SNN_I_point_subspace}
\begin{center}
\resizebox{\columnwidth}{!}{%
\small
\begin{tabular}{cccccccc}
\hline
Points & M & 20 & 40 & 60 & 80 & 100 & 300 \\
\hline
20 & $\|e\|_{ {L}^2}$ & 4.43e-01 & 1.38e-01 & 8.32e-02 & 2.20e-01 & 1.30e-01 & 7.30e-02 \\
& epochs & 683 & 713 & 964 & 837 & 954 & 1890 \\
30 & $\|e\|_{ {L}^2}$ & 4.89e-02 & 2.60e-01 & 2.21e-03 & 7.44e-04 & 2.89e-04 & 1.30e-03 \\
& epochs & 571 & 689 & 974 & 1187 & 1558 & 3595 \\
40 & $\|e\|_{ {L}^2}$ & 1.94e-02 & 1.27e-03 & 6.00e-05 & 9.16e-06 & 4.22e-05 & 1.38e-05 \\
& epochs & 432 & 731 & 1087 & 1244 & 1591 & 3246 \\
50 & $\|e\|_{ {L}^2}$ & 3.22e-02 & 5.52e-05 & 3.49e-05 & 6.87e-06 & 3.04e-06 & 3.53e-06 \\
& epochs & 407 & 686 & 1063 & 1421 & 1706 & 3484 \\
60 & $\|e\|_{ {L}^2}$ & 1.19e-02 & 1.05e-04 & 1.46e-06 & 6.59e-07 & 6.60e-07 & 2.54e-06 \\
& epochs & 402 & 688 & 897 & 1270 & 1735 & 3293 \\
80 & $\|e\|_{ {L}^2}$ & 2.94e-02 & 5.52e-05 & 4.91e-07 & 4.04e-07 & 5.93e-07 & 2.77e-07 \\
& epochs & 391 & 689 & 884 & 1210 & 1689 & 3776 \\
100 & $\|e\|_{ {L}^2}$ & 3.34e-02 & 5.96e-05 & 5.61e-07 & 4.77e-07 & 5.73e-07 & 8.41e-07 \\
& epochs & 400 & 718 & 885 & 1293 & 1672 & 3866 \\
300 & $\|e\|_{ {L}^2}$ & 3.97e-02 & 5.97e-05 & 5.37e-07 & 4.86e-07 & 5.77e-07 & 8.02e-07 \\
& epochs & 397 & 714 & 886 & 1301 & 1684 & 3708 \\
\hline
\end{tabular}
} 
\end{center}
\end{table}

\begin{figure}[htbp]
\centering
\subcaptionbox{}{\includegraphics[width=0.45\linewidth]{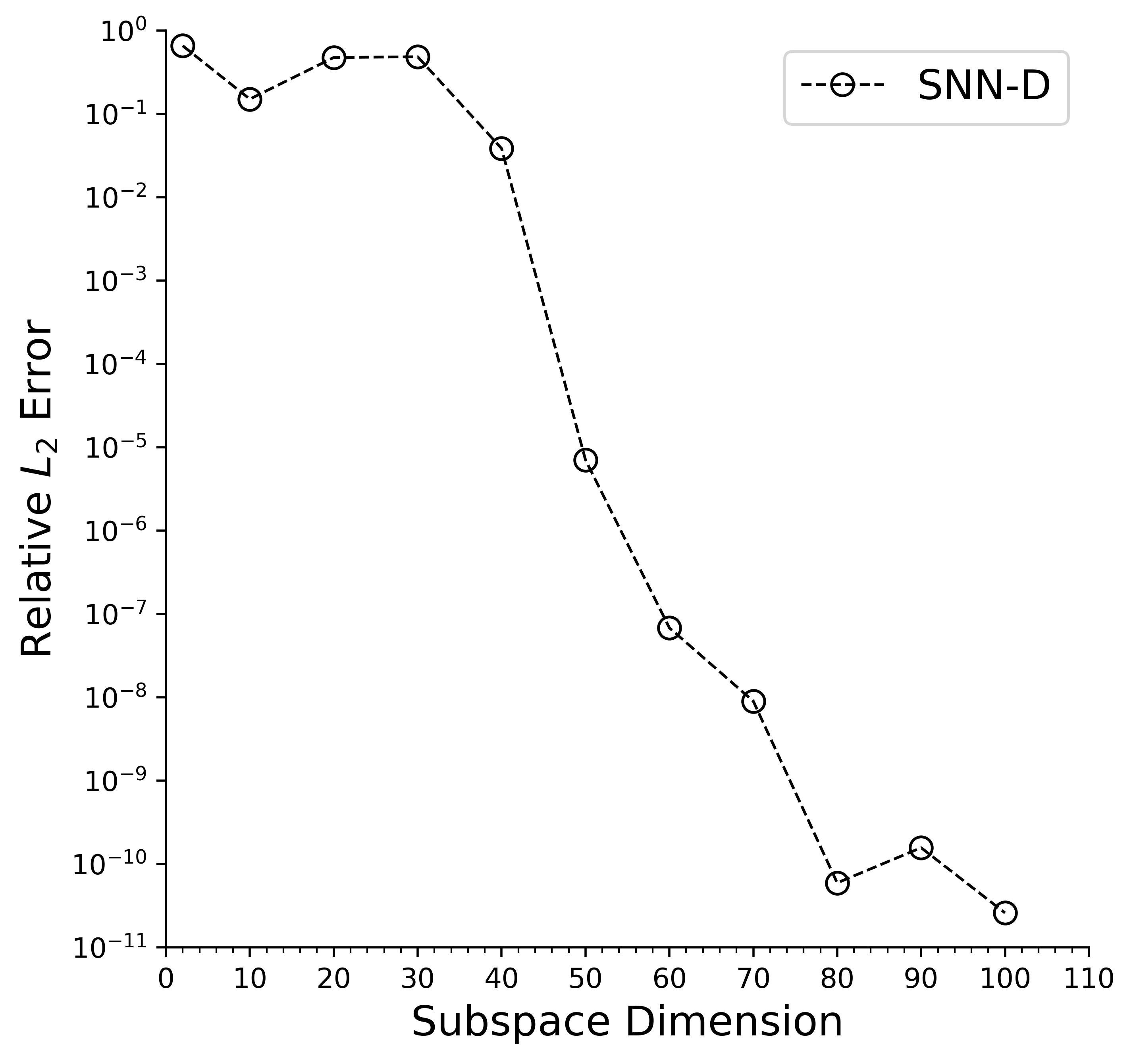}}
\hspace{0.3cm}  
\subcaptionbox{}{\includegraphics[width=0.45\linewidth]{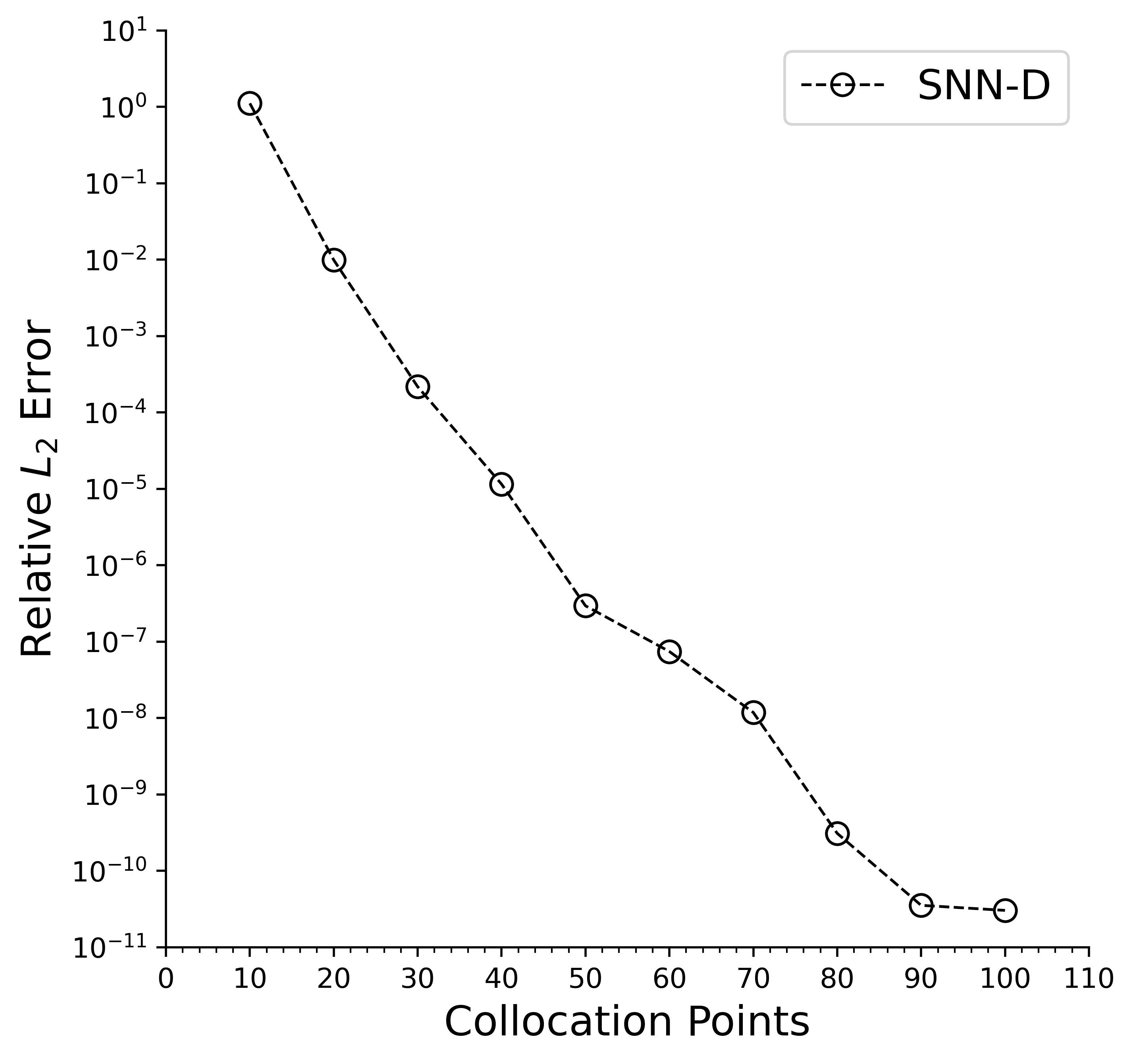}}
\caption{Error variation with subspace dimension at a fixed number of 1000 sampling points and error variation with the number of sampling points at a fixed subspace dimension of 300 for SNN-D.}
\label{1dhz_SNN_D_subspace_fig}
\end{figure}

\begin{figure}[htbp]
\centering
\subcaptionbox{}{\includegraphics[width=0.45\linewidth]{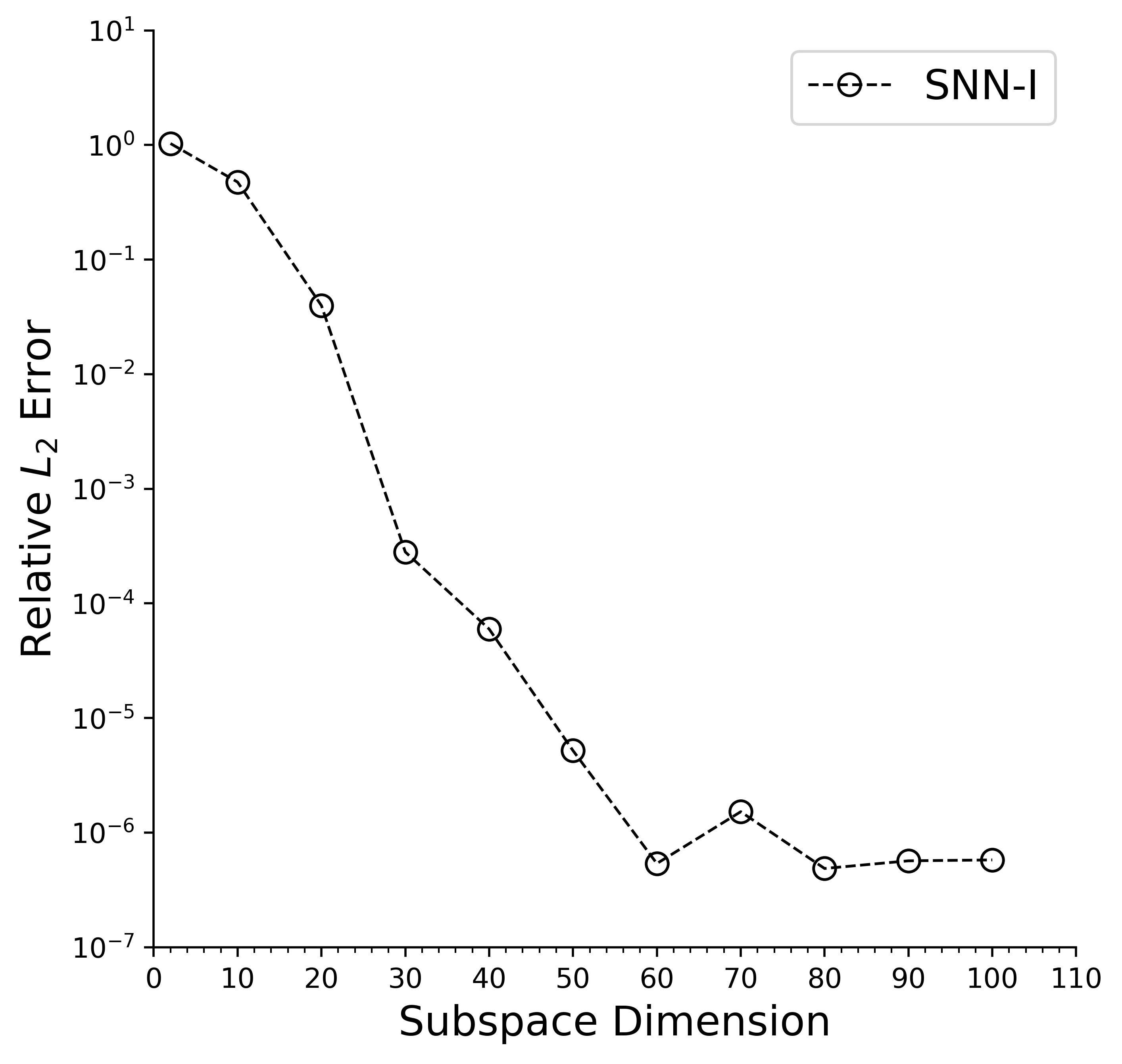}}
\hspace{0.3cm}  
\subcaptionbox{}{\includegraphics[width=0.45\linewidth]{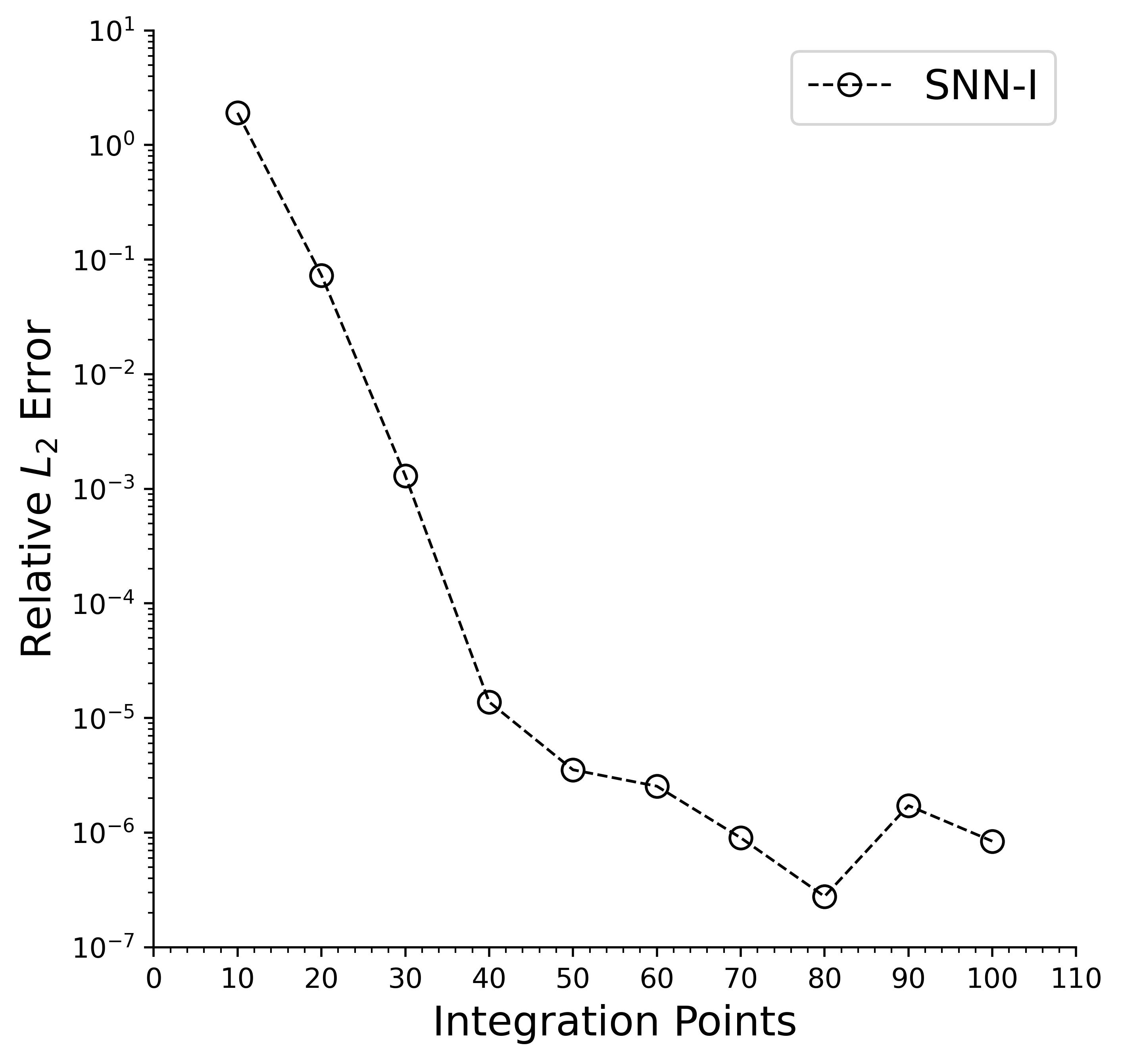}}
\caption{Error variation with subspace dimension at a fixed number of 300 integration points and error variation with the number of integration points at a fixed subspace dimension of 300 for SNN-I.}
\label{1dhz_SNN_I_subspace_fig}
\end{figure}

Now we  examine the variation pattern of errors. First, we fix the number of sampling points, and examine the error variation with subspace dimension. Then, we fix the subspace dimension, and examine the error variation with  the number of sampling points.
Figure \ref{1dhz_SNN_D_subspace_fig} illustrates the error variation with subspace dimension for  1000 sampling points, and the error variation with the number of  sampling points for a fixed subspace dimension of 300 for SNN-D. Similarly, Figure \ref{1dhz_SNN_I_subspace_fig} shows the variation pattern of errors for SNN-I.

Tables \ref{1dhz_SNN_D_netsize_subspace} and \ref{1dhz_SNN_I_netsize_subspace} present the results of SNN-D and SNN-I for different hidden layer depth and subspace dimension, maintaining a uniform width of 100 for each hidden layer. SNN-D utilizes a uniform sampling of 1000 points. SNN-I adopts the composite Gaussian quadrature rule, segmenting the domain into 30 subintervals with 10 points in each. 
When the number of hidden layer is small, even if the  subspace dimension is increased, SNN-D and SNN-I can only achieve the accuracy of $10^{-7}$ and $10^{-4}$, respectively. 
This indicates a constrained expressive capacity of the neural networks with few hidden layer, particularly when there is no hidden layer.
When there has three or four  hidden layers, the accuracy of method stabilizes, SNN-D achieves the accuracy of $10^{-11}$ and SNN-I achieves the accuracy of  $10^{-7}$. 
Furthermore, a trend can be  observed that the number of training epochs significantly decreases when the number of hidden layer increases. 
When there has no hidden layer, the number of training epochs reaches the maximum limit, and the accuracy of SNN-I can not be improved by increasing the subspace dimension. This indicates that the hidden layers are necessary in order to improve the accuracy and robustness. 

\begin{table}[!htbp]
\caption{The errors and epochs for SNN-D on different hidden layer depth.}\label{1dhz_SNN_D_netsize_subspace}
\begin{center}
\small
\begin{tabular}{ccccc}
\hline
Hidden Layer & M  & 60 & 100 & 300 \\
\hline
0 & $\|e\|_{ {L}^2}$ & 3.87e-04 & 1.12e-05 & 2.13e-07 \\
& epochs & 5000 & 5000 & 5000 \\
1 & $\|e\|_{ {L}^2}$ & 3.29e-09 & 1.70e-08 & 4.21e-08 \\
& epochs & 1182 & 1103 & 1791 \\
2 & $\|e\|_{ {L}^2}$ & 8.93e-08 & 1.07e-10 & 1.65e-10 \\
& epochs & 594 & 541 & 1180 \\
3 & $\|e\|_{ {L}^2}$ & 9.97e-08 & 1.33e-10 & 1.25e-10 \\
& epochs & 713 & 659 & 2523 \\
4 & $\|e\|_{ {L}^2}$ & 6.79e-08 & 2.58e-11 & 1.79e-11 \\
& epochs & 729 & 1184 & 1960 \\
5 & $\|e\|_{ {L}^2}$ & 1.46e-07 & 2.20e-10 & 4.99e-11 \\
& epochs & 726 & 943 & 978 \\
6 & $\|e\|_{ {L}^2}$ & 4.91e-07 & 1.15e-10 & 2.53e-10 \\
& epochs & 619 & 773 & 768 \\
7 & $\|e\|_{ {L}^2}$ & 2.47e-07 & 8.51e-11 & 3.45e-11 \\
& epochs & 588 & 510 & 606 \\
8 & $\|e\|_{ {L}^2}$ & 1.80e-05 & 4.99e-11 & 7.59e-11 \\
& epochs & 552 & 541 & 541 \\
\hline
\end{tabular}
\end{center}
\end{table}

\begin{table}[!htbp]
\caption{The errors and epochs for SNN-I on different hidden layer depth.}\label{1dhz_SNN_I_netsize_subspace}
\begin{center}
\small
\begin{tabular}{ccccc}
\hline
Hidden Layer & M & 60 & 100 & 300 \\
\hline
0 & $\|e\|_{ {L}^2}$ & 7.91e-04 & 3.81e-03 & 5.52e-04 \\
& epochs & 5000 & 5000 & 5000 \\
1 & $\|e\|_{ {L}^2}$ & 1.94e-04 & 1.81e-04 & 2.95e-05 \\
& epochs & 1169 & 1473 & 1694 \\
2 & $\|e\|_{ {L}^2}$ & 3.12e-07 & 3.34e-07 & 4.32e-07 \\
& epochs & 598 & 631 & 1245 \\
3 & $\|e\|_{ {L}^2}$ & 4.34e-07 & 1.39e-06 & 6.51e-07 \\
& epochs & 647 & 797 & 4173 \\
4 & $\|e\|_{ {L}^2}$ & 5.37e-07 & 5.77e-07 & 8.02e-07 \\
& epochs & 886 & 1684 & 3708 \\
5 & $\|e\|_{ {L}^2}$ & 7.54e-07 & 5.08e-07 & 1.43e-06 \\
& epochs & 853 & 831 & 997 \\
6 & $\|e\|_{ {L}^2}$ & 1.80e-06 & 1.18e-06 & 1.26e-06 \\
& epochs & 688 & 837 & 1143 \\
7 & $\|e\|_{ {L}^2}$ & 4.13e-06 & 4.01e-07 & 2.36e-07 \\
& epochs & 633 & 636 & 694 \\
8 & $\|e\|_{ {L}^2}$ & 3.10e-06 & 1.50e-06 & 4.96e-07 \\
& epochs & 501 & 580 & 507 \\
\hline
\end{tabular}
\end{center}
\end{table}

\subsection{Poisson equation}
\label{Poisson equation}
In the second test, we consider the two-dimensional  Poisson equation  on the domain $\Omega=(0,1) \times(0,1)$,
\begin{eqnarray}
    \label{poisson}
    -\Delta u=f(x, y), & (x, y) \in \Omega, \\ u(x, y)=g(x, y), & (x, y) \in \partial \Omega.
\end{eqnarray}
We choose  $f (x,y)$ and $g (x,y)$ such that Eq. (\ref{poisson}) has the following solution,
\begin{eqnarray*}
u(x) = \sin\left(\pi x\right) \sin\left(\pi y \right).
\end{eqnarray*}

For SNN-D, \(32 \times 32\) points are uniformly sampled throughout the domain, ensuring an uniform distribution of 32 points along each boundary. For SNN-I, a two-dimensional composite Gaussian quadrature formula is used, segmenting each dimension into 8 subintervals and allocating 4 points to each. This strategy leads to 1024 sampling points in the interior and 128 sampling points on the boundary.

\begin{figure}[htbp]
\centering
\subcaptionbox{SNN-D}{\includegraphics[width=0.45\linewidth]{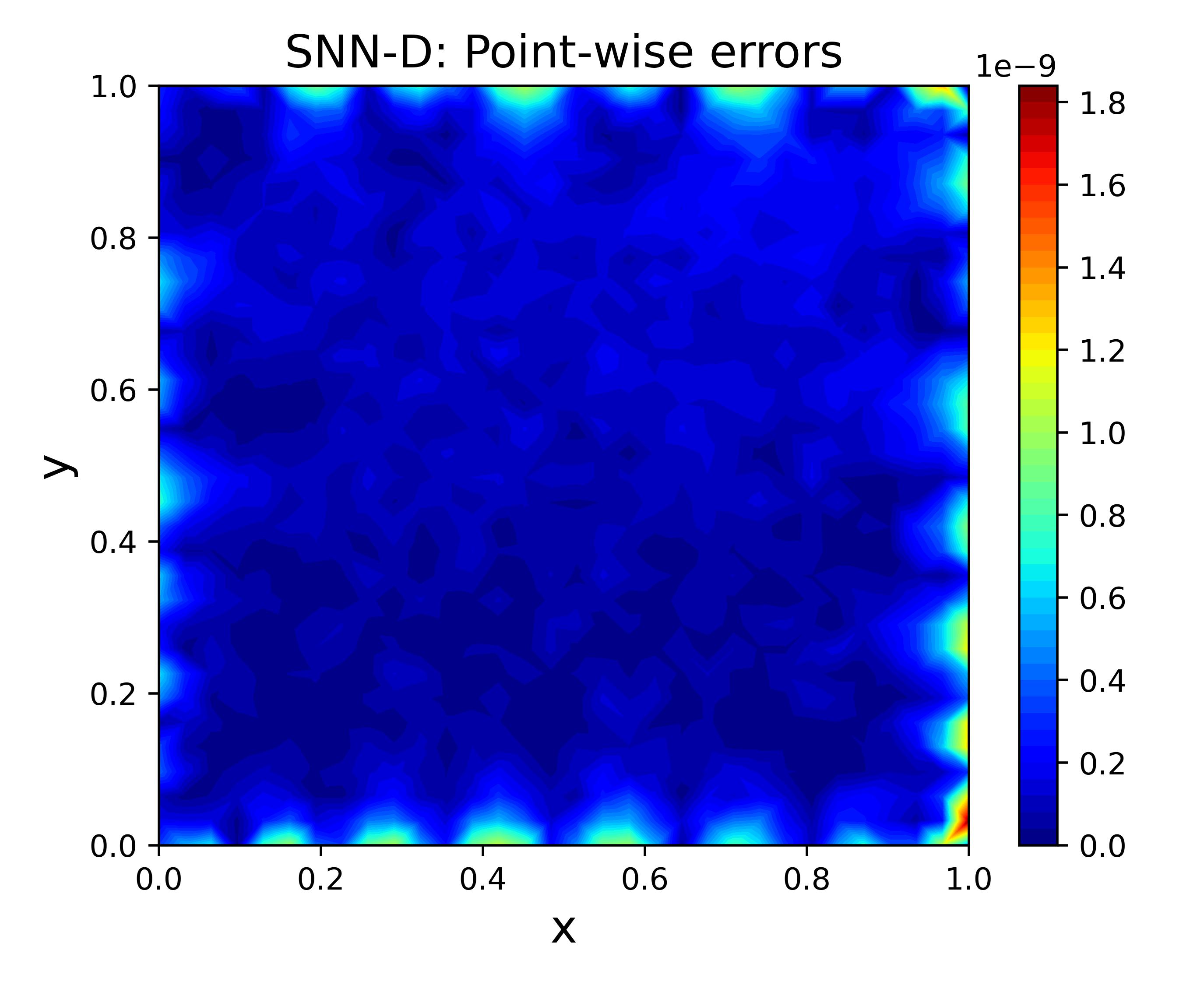}}
\hspace{0.3cm}  
\subcaptionbox{SNN-I}{\includegraphics[width=0.45\linewidth]{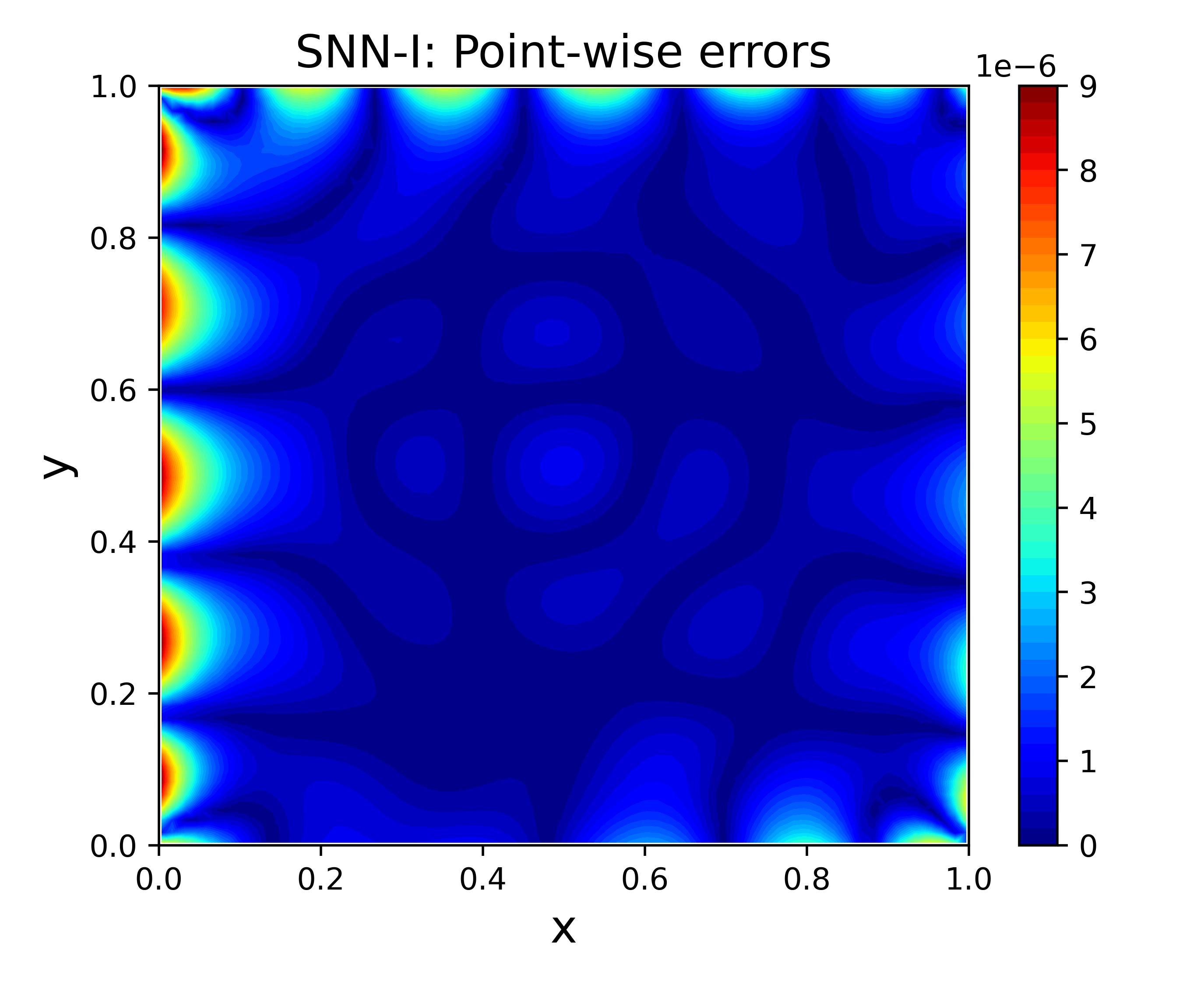}}
\caption{  Point-wise errors   of SNN-D and SNN-I  for Poisson equation.}
\label{2d_ps_point_error}
\end{figure}

Figure \ref{2d_ps_point_error} illustrates the point-wise errors  using SNN-D and SNN-I for  Poisson equation.  Table \ref{2d_ps_error} presents the relative \(L^2\) errors for various methods including SNN-D, SNN-I, PINN, DGM and ELM. This example serves to demonstrate the adaptivity of our algorithms. Notably, after 50000 training epochs, the relative \(L^2\) errors of PINN and DGM are 2.71e-04 and 1.81e-03, respectively. However, SNN-D and SNN-I achieve relative \(L_2\) errors of 5.37e-10 and 2.97e-06, respectively, with significantly fewer epochs, specifically 276 and 275. Compared to 1D Helmholtz equation, 2D Poisson equation requires fewer training epochs, typically ranging from 40 to 500. The relative \(L_2\) error of ELM is 6.10e-11, suggesting that the subspace initially provides a fairly accurate approximation to the solution space. This highlights the adaptability of our algorithm. In this test, these initial parameters can approximate the solution space  well, which leads to a rapid decrease of the loss function. Conversely, in Section \ref {Helmholtz equation}, due to the initial parameters poorly approximating the solution space, the loss function decreases slowly, which leads to more training epochs. Thus, our approach adaptively adjusts the number of training epochs required.

\begin{table}[!htbp]
\caption{The errors and epochs of different methods for Poisson equation.}\label{2d_ps_error}
\begin{center}
\small
\begin{tabular}{cccccc}\hline
Method & $\|e\|_{L^2}$ & $\|e\|_{L^\infty}$ & epochs\\
\hline
PINN & 2.71e-04 & 7.33e-03 & 50000 \\
ELM & 6.10e-11 & 1.60e-10 & 0 \\
SNN-D & 5.37e-10 & 1.81e-09 & 276 \\
DGM & 1.81e-03 & 3.26e-03 & 50000 \\
SNN-I & 2.97e-06 & 7.97e-06 & 275 \\
\hline
\end{tabular}
\end{center}
\end{table}

Tables \ref{2dps_SNN_D_point_subspace} and \ref{2dps_SNN_I_point_subspace} present the performance of SNN-D and SNN-I across various numbers of sampling points and subspace dimension, maintaining a networks architecture of four hidden layers with 100 neurons each. SNN-D utilizes uniform sampling, while SNN-I adopts the composite Gaussian quadrature, allocating 4 points per subinterval in each direction and increasing subinterval numbers for enhanced precision. Each method preserves equal sampling points on boundary as in a single direction. 
Similar to the 1D Helmholtz equation, 
it becomes apparent that the error decreases with increasing both sampling points and subspace dimension. 
With a subspace dimension of 50, the error primarily ranges between \(10^{-2}\) and \(10^{-3}\).
Expanding the subspace dimension to 100, SNN-D can achieve the accuracy of  \(10^{-4}\) and SNN-I can achieve the accuracy of $10^{-6}$. With enough sampling points and further expansion of the subspace dimension,  SNN-D achieves the accuracy of $10^{-10}$, and the accuracy of SNN-I stables at \(10^{-6}\).

\begin{table}[!htbp]
\caption{The errors and epochs for SNN-D across various numbers of sampling points and subspace dimension $M$ for Poisson equation.} 
\label{2dps_SNN_D_point_subspace}
\begin{center}
\resizebox{\columnwidth}{!}{%
\small
\begin{tabular}{ccccccccc}
\hline
Points & $M$ & 50 & 100 & 150 & 200 & 250 & 300 & 500 \\
\hline
$8\times8$ & $\|e\|_{ {L}^2}$ & 7.24e-03 & 7.06e-05 & 1.82e-05 & 2.75e-05 & 2.77e-05 & 1.91e-05 & 2.55e-05 \\
 & epochs & 206 & 224 & 249 & 243 & 235 & 245 & 315 \\
$12\times12$ & $\|e\|_{ {L}^2}$ & 7.06e-03 & 1.27e-04 & 5.95e-07 & 9.28e-08 & 8.56e-08 & 9.59e-08 & 7.06e-08 \\
 & epochs & 217 & 233 & 251 & 252 & 243 & 263 & 344 \\
$16\times16$ & $\|e\|_{ {L}^2}$ & 6.86e-03 & 1.23e-04 & 1.68e-06 & 2.93e-08 & 4.26e-09 & 1.91e-09 & 6.85e-10 \\
 & epochs & 223 & 236 & 251 & 258 & 246 & 271 & 354 \\
$24\times24$ & $\|e\|_{ {L}^2}$ & 8.69e-03 & 1.06e-04 & 1.82e-06 & 4.19e-08 & 3.00e-09 & 5.24e-10 & 2.66e-10 \\
 & epochs & 229 & 239 & 250 & 262 & 246 & 275 & 376 \\
$32\times32$ & $\|e\|_{ {L}^2}$ & 1.05e-02 & 1.07e-04 & 1.85e-06 & 3.07e-08 & 4.28e-09 & 5.37e-10 & 2.83e-10 \\
 & epochs & 233 & 240 & 249 & 264 & 246 & 276 & 396 \\
$48\times48$ & $\|e\|_{ {L}^2}$ & 1.32e-02 & 1.26e-04 & 1.97e-06 & 2.13e-08 & 4.71e-09 & 1.08e-09 & 2.92e-10 \\
 & epochs & 236 & 241 & 247 & 264 & 244 & 277 & 404 \\
\hline
\end{tabular}
} 
\end{center}
\end{table}

\begin{table}[!htbp]
\caption{The errors and epochs for SNN-I across various numbers of sampling points and subspace dimension $M$ for Poisson equation.} 
\label{2dps_SNN_I_point_subspace}
\begin{center}
\resizebox{\columnwidth}{!}{%
\small
\begin{tabular}{ccccccccc}
\hline
Points & $M$ & 50 & 100 & 150 & 200 & 250 & 300 & 500 \\
\hline
$4\times4$ & $\|e\|_{ {L}^2}$ & 1.70e-02 & 3.87e-03 & 1.2e-03 & 3.15e-03 & 4.32e-03 & 1.59e-03 & 2.19e-03 \\
 & epochs & 48 & 53 & 54 & 49 & 48 & 60 & 59 \\
$8\times8$ & $\|e\|_{ {L}^2}$ & 1.08e-03 & 6.85e-06 & 7.13e-06 & 1.14e-05 & 1.21e-05 & 8.96e-06 & 1.42e-05 \\
 & epochs & 243 & 242 & 243 & 263 & 240 & 275 & 417 \\
$12\times12$ & $\|e\|_{ {L}^2}$ & 1.66e-03 & 3.84e-06 & 2.34e-06 & 2.33e-06 & 1.67e-06 & 4.45e-06 & 2.54e-06 \\
 & epochs & 242 & 241 & 242 & 262 & 239 & 275 & 410 \\
$16\times16$ & $\|e\|_{ {L}^2}$ & 1.59e-03 & 3.18e-06 & 1.57e-06 & 2.22e-06 & 1.79e-06 & 2.72e-06 & 2.63e-06 \\
 & epochs & 242 & 241 & 242 & 262 & 239 & 275 & 411 \\
$24\times24$ & $\|e\|_{ {L}^2}$ & 1.60e-03 & 3.75e-06 & 2.75e-06 & 1.51e-06 & 1.91e-06 & 3.02e-06 & 1.97e-06 \\
 & epochs & 242 & 241 & 242 & 262 & 239 & 275 & 411 \\
$32\times32$ & $\|e\|_{ {L}^2}$ & 1.55e-03 & 4.00e-06 & 2.76e-06 & 1.51e-06 & 1.87e-06 & 2.97e-06 & 1.95e-06 \\
 & epochs & 242 & 241 & 242 & 262 & 239 & 275 & 411 \\
\hline
\end{tabular}
} 
\end{center}
\end{table}

\begin{figure}[htbp]
\centering
\subcaptionbox{}{\includegraphics[width=0.45\linewidth]{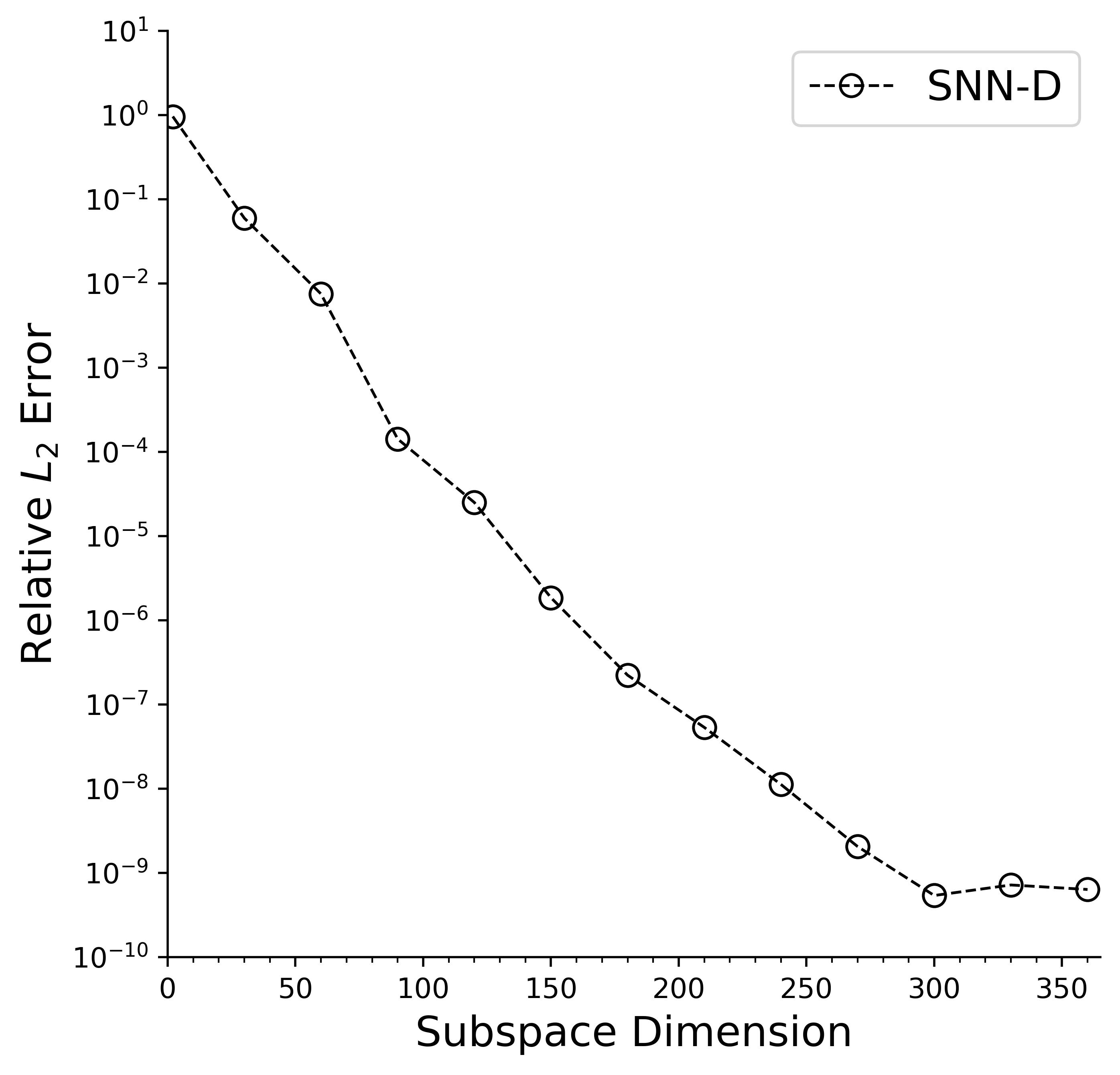}}
\hspace{0.3cm}  
\subcaptionbox{}{\includegraphics[width=0.45\linewidth]{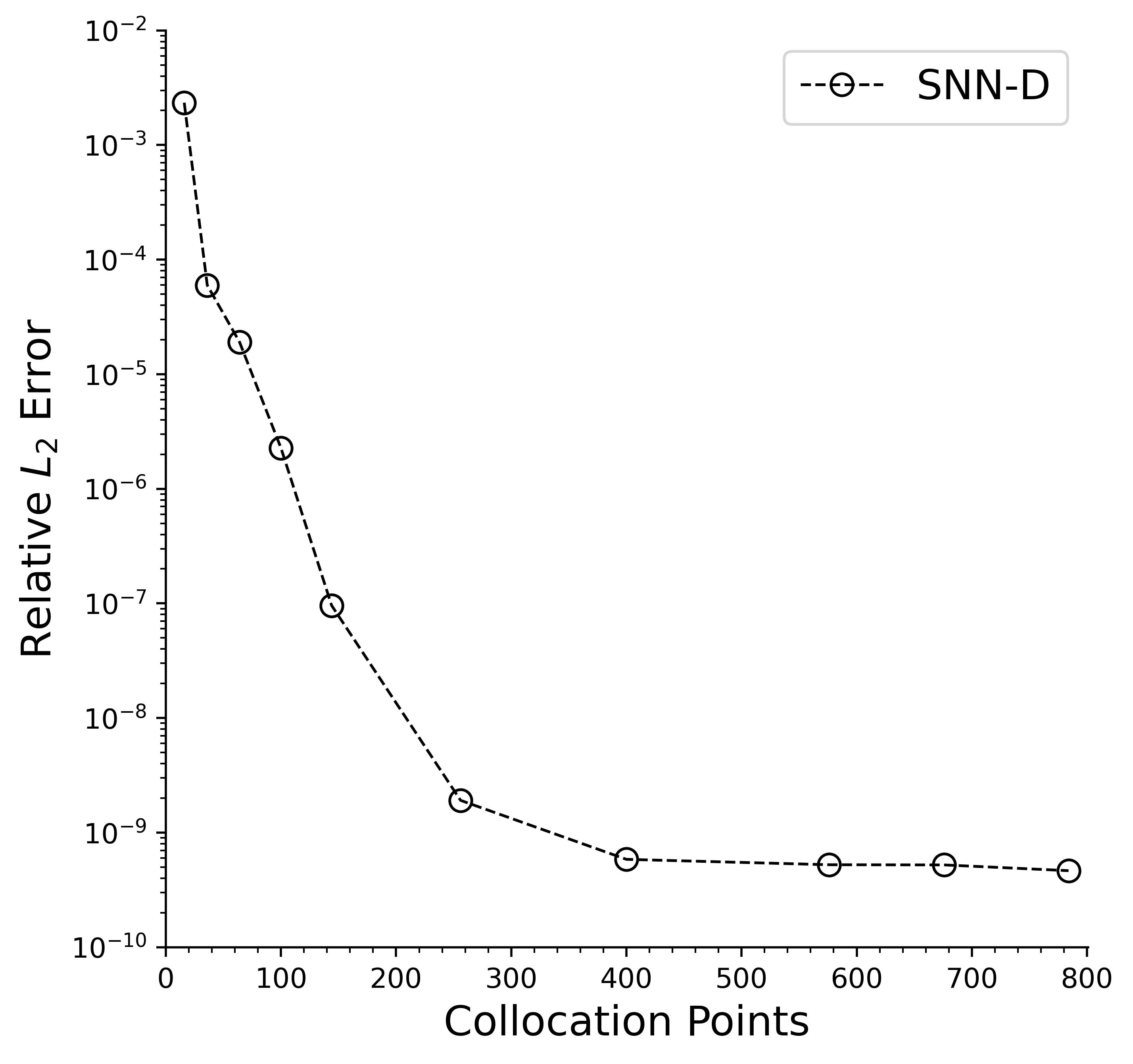}}
\caption{ Error variation with subspace dimension at a fixed number of \(32 \times 32\) sampling points and error variation with the number of sampling points at a fixed subspace dimension of 300 for SNN-D.}
\label{2dps_SNN_D_subspace_fig}
\end{figure}

\begin{figure}[htbp]
\centering
\subcaptionbox{}{\includegraphics[width=0.45\linewidth]{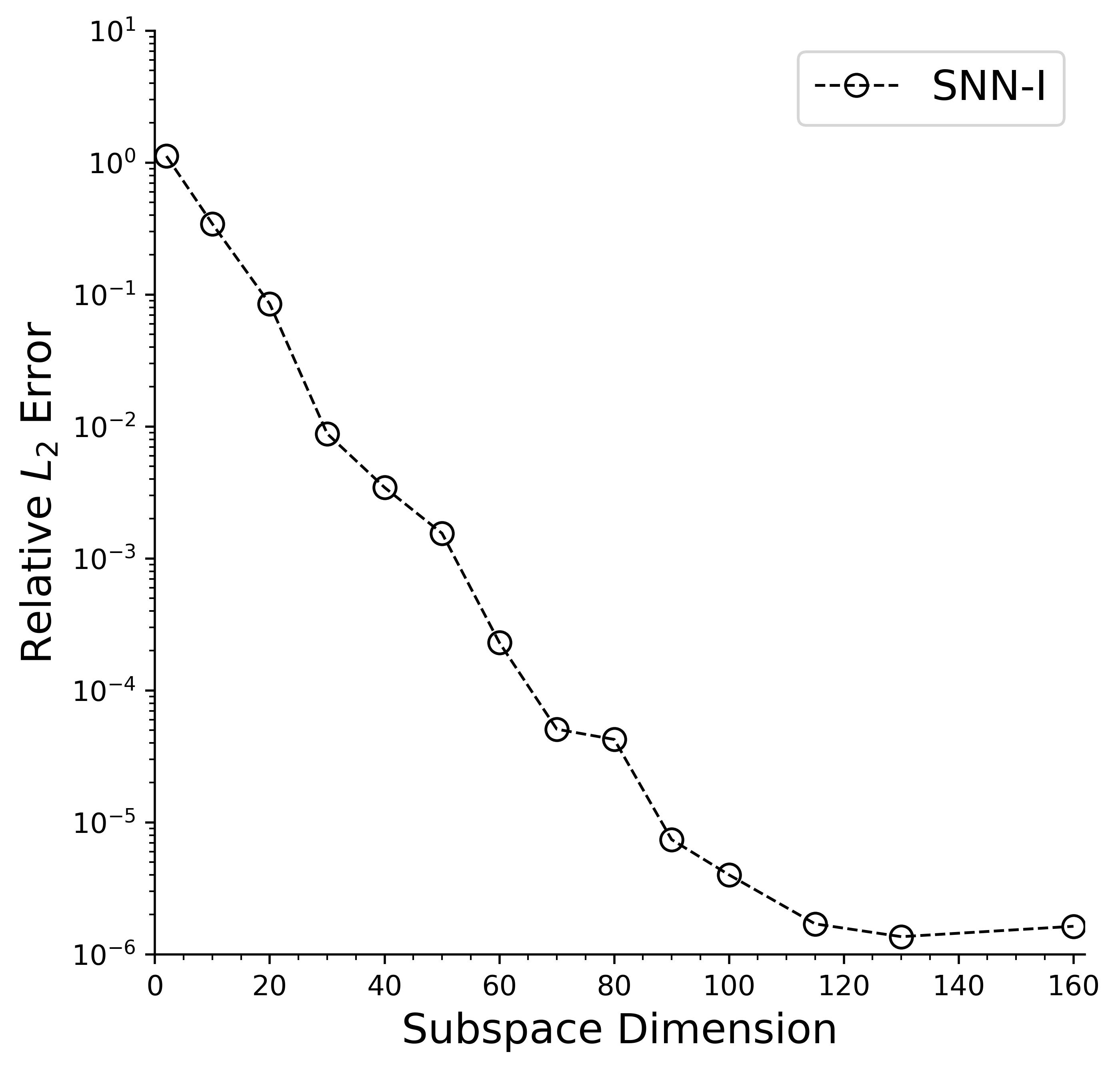}}
\hspace{0.3cm}  
\subcaptionbox{}{\includegraphics[width=0.45\linewidth]{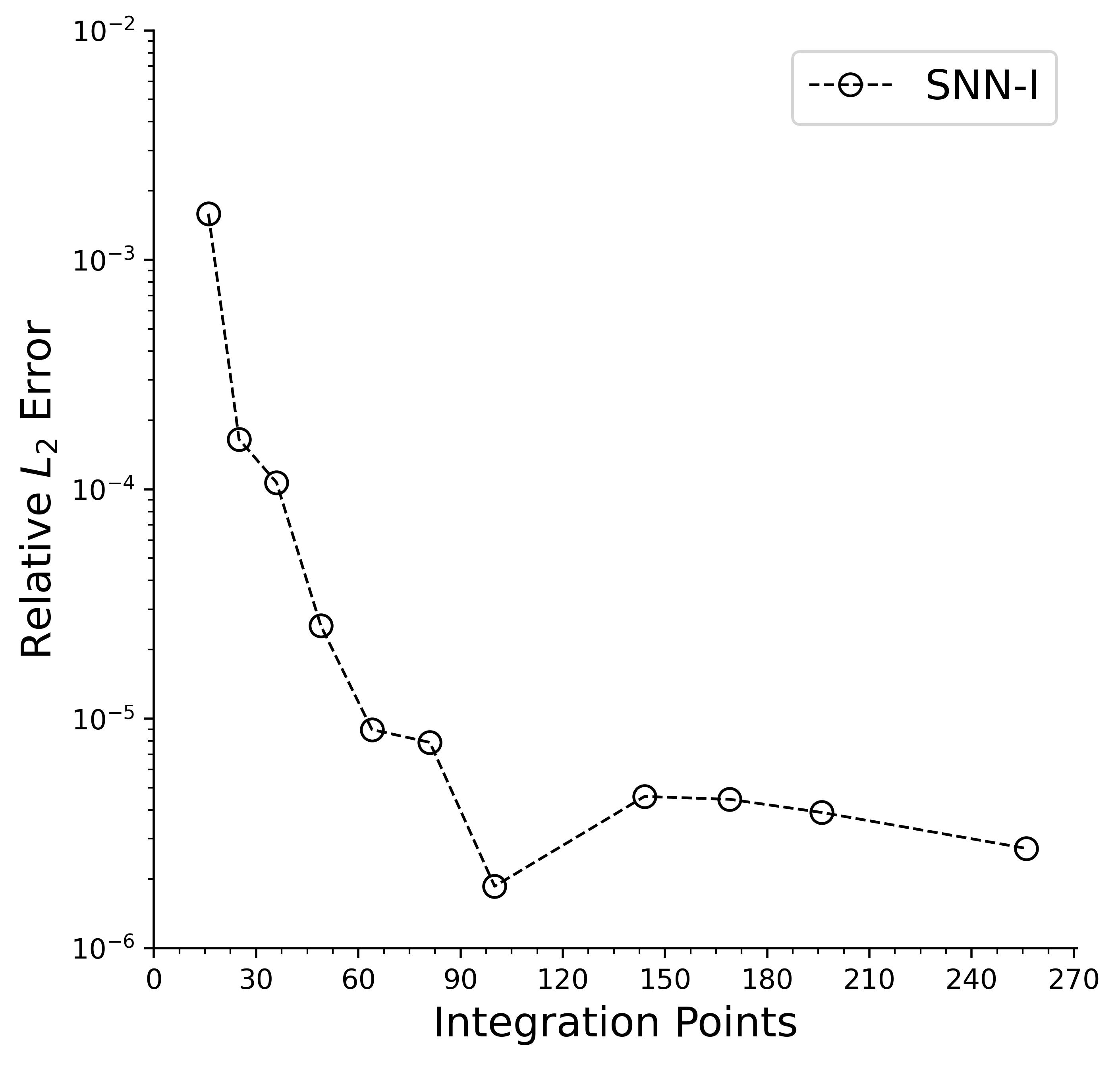}}
\caption{Error variation with subspace dimension at a fixed number of \(32 \times 32\) integration points and error variation with the number of integration points at a fixed subspace dimension of 300 for SNN-I.}
\label{2dps_SNN_I_subspace_fig}
\end{figure}

Figure \ref{2dps_SNN_D_subspace_fig} illustrates the error variation with subspace dimension for \(32 \times 32\) sampling points, and the error variation with the number of sampling points for a fixed subspace dimension of 300 for SNN-D. 
Figure \ref{2dps_SNN_I_subspace_fig} illustrates the error variation with subspace dimension for \(32 \times 32\) Gaussian integration points, and  the error variation with the number of Gaussian integration points for a fixed subspace dimension of 300 for SNN-I.
The number of training epochs for these experiments is from  48 to 417, significantly lower than that of PINN and DGM.

\subsection{Advection equation}
\label{Advection equation}
Next, we test our  method for the advection equation. 
Consider  the following advection equation in space-time domain $\Omega =(a,b)\times(0,T)$, 
\begin{equation}
 \frac{\partial u}{\partial t}-c \frac{\partial u}{\partial x}=0.
\end{equation}
The initial condition $u(x, 0)=h(x)$, and we impose periodic boundary condition $u\left(a, t\right)=u\left(b, t\right)$.  In this example, we take $a=-1$, $b=1$, $T=1$, $c=-2$, $u(x,t)=sin(\pi(x-2t))$.

For the deep learning methods such as PINN and DGM, the advection equation falls into a category where training progress tends to be slow. This is due to the advection equation's strong reliance on initial boundary conditions, which causes the imbalance between PDE loss and initial boundary loss in the training process. 
Some weighted techniques have been used to correct this imbalance. However, how to determine these weights is a challenging problem. We do not use the initial boundary condition in the loss function, hence we do not
need to introduce these weights.
Through this example, we will show the accuracy and effect of our method without the initial boundary loss for solving time-dependent equation.

For SNN-D, \(100 \times 100\) points are uniformly sampled across the domain, ensuring 100 points are uniformly distributed along each boundary and 500 points at the initial time. For SNN-I, a two-dimensional composite Gaussian quadrature formula is applied, dividing each dimension into 10 subintervals with 10 points each. The initial condition is divided into 50 subintervals. This approach results in 10000 sampling points in the interior, 100 points on each periodic boundary, and 500 points at the initial time.

\begin{figure}[htbp]
\centering
\subcaptionbox{SNN-D}{\includegraphics[width=0.45\linewidth]{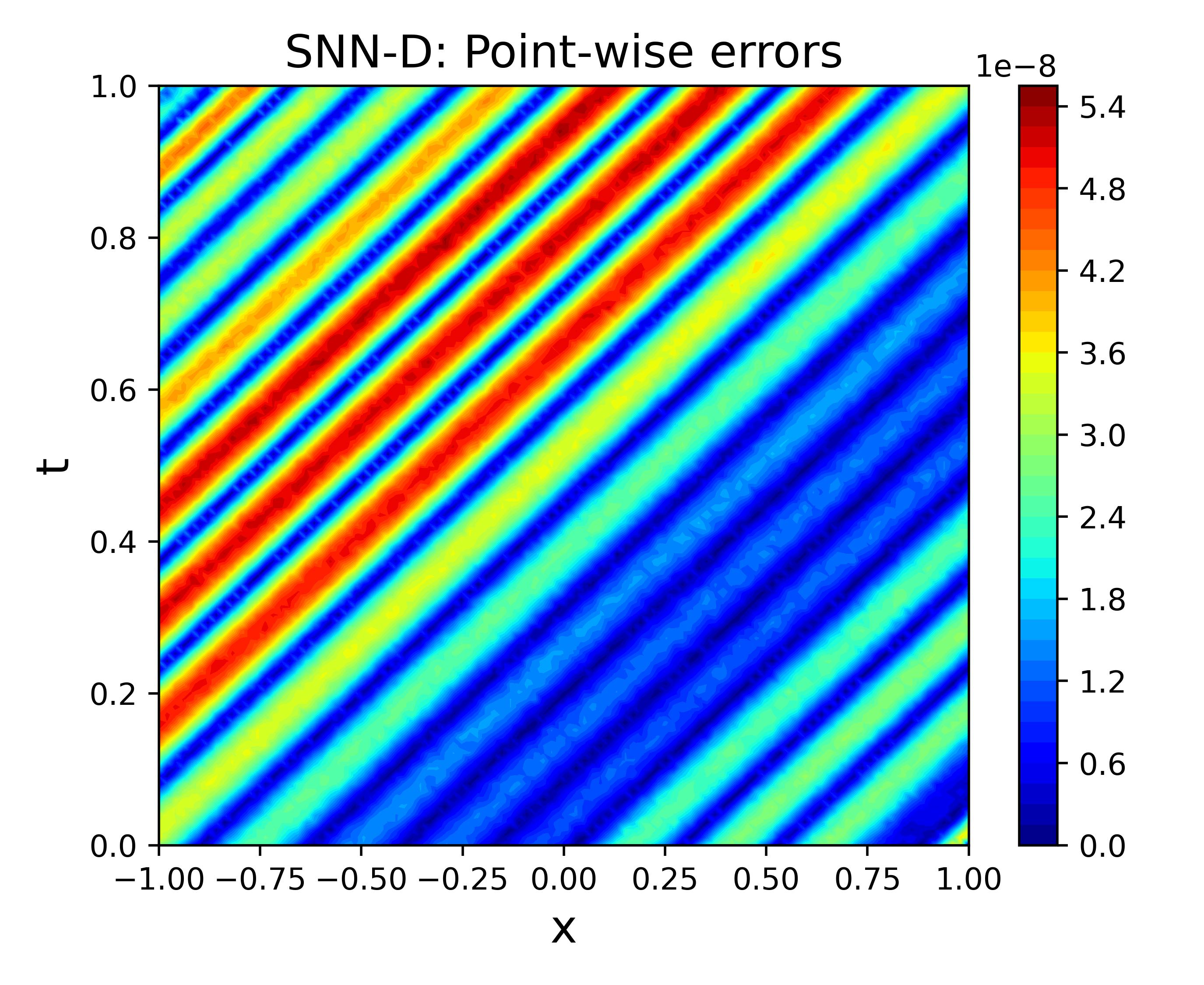}}
\hspace{0.3cm}  
\subcaptionbox{SNN-I}{\includegraphics[width=0.45\linewidth]{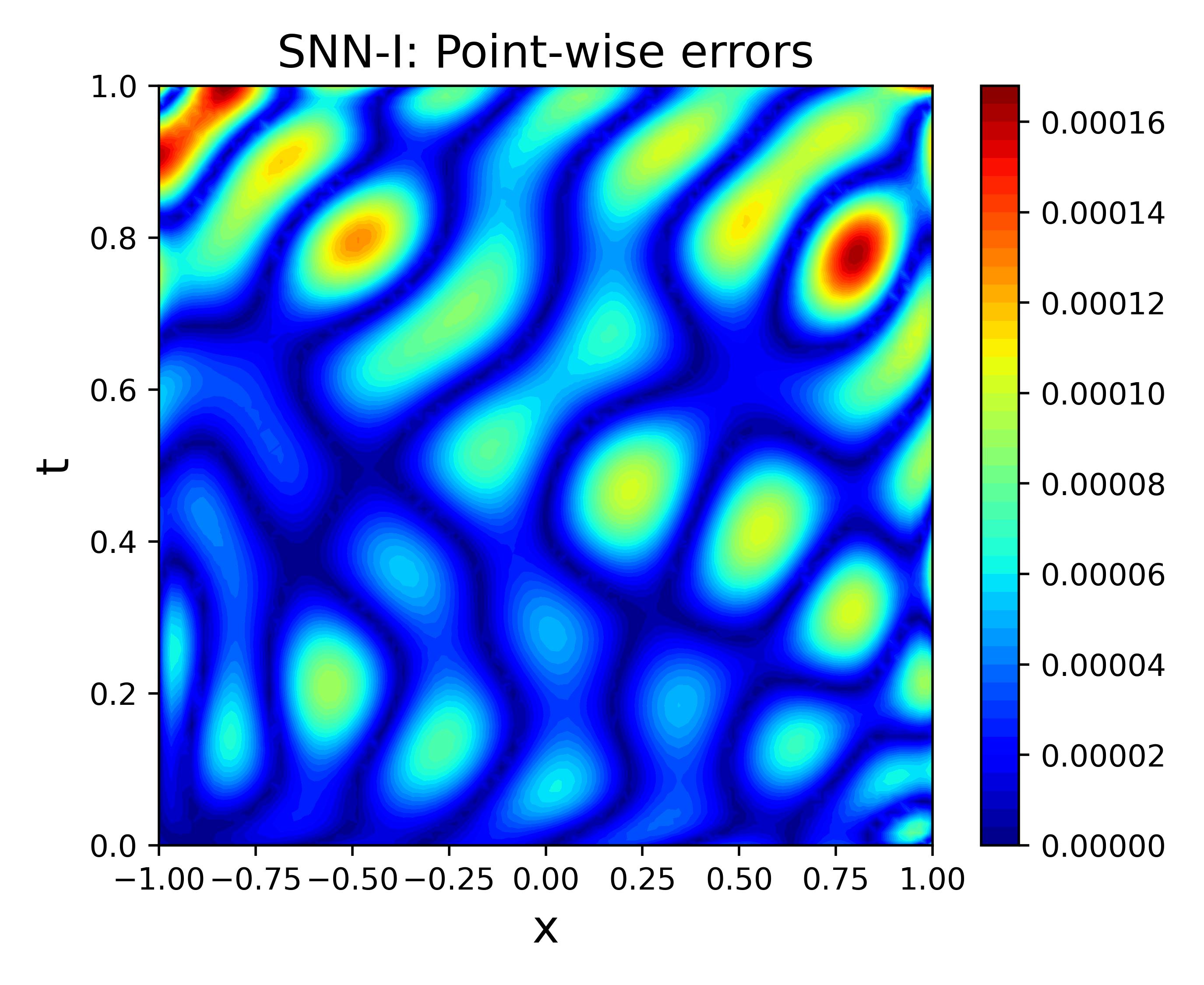}}
\caption{  Point-wise errors  of SNN-D and SNN-I for the advection equation.}
\label{ad_point_error}
\end{figure}

Figure \ref{ad_point_error} illustrates the point-wise errors for both SNN-D and SNN-I in this example. Table \ref{ad_error} presents the relative \(L_2\) errors for  two methods. It is observed that the errors for SNN-D and SNN-I are 3.62e-07 and 6.78e-05, respectively. This indicates that our methods continue to maintain precision and efficiency for advection equation.

\begin{table}[!htbp]
\caption{The errors and epochs of SNN-D and SNN-I for the advection equation.}\label{ad_error}
\begin{center}
\small
\begin{tabular}{cccccc}\hline
  Method & $\|e\|_{ {L}^2}$ & $\|e\|_{L^{\infty}}$ & epochs\\
\hline
SNN-D &   3.62e-08    &   5.48e-08   &    84   \\
SNN-I &   6.78e-05    &   1.62e-04   &    84   \\
\hline
\end{tabular}
\end{center}
\end{table}

Tables \ref{ad_D_withoutBC} and \ref{ad_D_withBC} give the errors and epochs of SNN-D without and with the initial boundary loss term, respectively. SNN-D without the initial boundary loss term  only need 73 to 97 epochs, whereas SNN-D with the initial boundary loss term need  3205 to 5000 epochs. This indicates that the training time differs by about 50 times between two approaches. When the subspace dimension is between 100 and 300, the error of SNN-D without the initial boundary loss term  generally falls below that of SNN-D with the initial boundary loss term, though the latter shows slightly higher accuracy when subspace dimension is 50 and 500.

\begin{table}[!htbp]
\caption{The errors and epochs of SNN-D without the initial boundary loss term for the advection equation.}\label{ad_D_withoutBC}
\begin{center}
\small
\begin{tabular}{cccccc}\hline
M & $50$ & $100$ & $200$ & $300$ & $500$\\
\hline
$\|e\|_{ {L}^2}$ & 2.44e-01 & 9.06e-04 & 2.18e-07 & 3.62e-08 & 3.85e-08 \\
epochs & 97 & 73 & 74 & 84 & 87 \\
\hline
\end{tabular}
\end{center}
\end{table}

\begin{table}[!htbp]
\caption{The errors and epochs of SNN-D with the initial boundary loss term for the advection equation.}\label{ad_D_withBC}
\begin{center}
\small
\begin{tabular}{cccccc}\hline
M & $50$ & $100$ & $200$ & $300$ & $500$\\
\hline
$\|e\|_{ {L}^2}$ & 2.44e-02 & 1.31e-03 & 5.38e-06 & 1.03e-07 & 3.17e-09 \\
epochs & 3205 & 3294 & 5000 & 5000 & 5000 \\
\hline
\end{tabular}
\end{center}
\end{table}

Tables \ref{ad_I_withoutBC} and \ref{ad_I_withBC} give the errors and epochs of SNN-I without and with the initial boundary loss term, respectively. we can see that the training time of SNN-I without the initial boundary loss term is significantly less than that of SNN-I with the initial boundary loss term.

\begin{table}[!htbp]
\caption{The errors and epochs of SNN-I without the initial boundary loss term for the advection equation.}\label{ad_I_withoutBC}
\begin{center}
\small
\begin{tabular}{cccccc}\hline
M & $50$ & $100$ & $200$ & $300$ & $500$\\
\hline
$\|e\|_{ {L}^2}$ & 1.00e-02 & 2.01e-04 & 1.10e-04 & 6.78e-05 & 1.51e-04 \\
epochs & 97 & 73 & 74 & 84 & 75 \\
\hline
\end{tabular}
\end{center}
\end{table}

\begin{table}[!htbp]
\caption{The errors and epochs of SNN-I with the initial boundary loss term for the advection equation.}\label{ad_I_withBC}
\begin{center}
\small
\begin{tabular}{cccccc}\hline
M & $50$ & $100$ & $200$ & $300$ & $500$\\
\hline
$\|e\|_{ {L}^2}$ & 3.81e-03 & 2.28e-04 & 1.32e-05 & 2.13e-05 & 2.80e-05 \\
epochs & 2391 & 2751 & 5000 & 3606 & 5000 \\
\hline
\end{tabular}
\end{center}
\end{table}

It is a challenging problem to balance the  PDE loss and initial boundary loss in the training process.
SNN-D and SNN-I with the initial boundary loss term suffer from the imbalance between the  PDE loss and initial boundary loss, which lead to a significant increase of training epochs, even the epochs reach the preset maximum number \(N_{max} = 5000\). However, SNN-D and SNN-I without the initial boundary loss term overcome this difficulty, and can complete the training with fewer epochs.

\subsection{Parabolic equation}
\label{Parabolic equation} 
Consider the following parabolic equation on the space-time domain $\Omega =(a,b)\times (0,T)$, 
\begin{eqnarray*}
 \frac{\partial u}{\partial t} - \frac{\partial^2 u}{\partial x^2} &=& f(x), \\
 u(a, t) &=& g_1(t), \\
 u(b, t) &=& g_2(t), \\
  u(x, 0) &=& h(x),
\end{eqnarray*}
where \(h(x)\) is the initial condition,  \(g_1(t)\) and \(g_2(t)\) are given functions.  We take
$a = 0$, $b = 1$, $T = 1$, $u(x,t) = 2e^{-t}\sin(\pi x)$.

For SNN-D, we uniformly sample \(50 \times 50\) points across the domain, with 50 points sampled uniformly at the initial time and each boundary. For SNN-I, we employ a two-dimensional composite Gaussian quadrature formula, dividing each dimension into 10 subintervals, with 5 points per subinterval. This strategy leads to 2500 sampling points in the interior,  and 150 sampling points on the boundary.

\begin{figure}[htbp]
\centering
\subcaptionbox{SNN-D}{\includegraphics[width=0.45\linewidth]{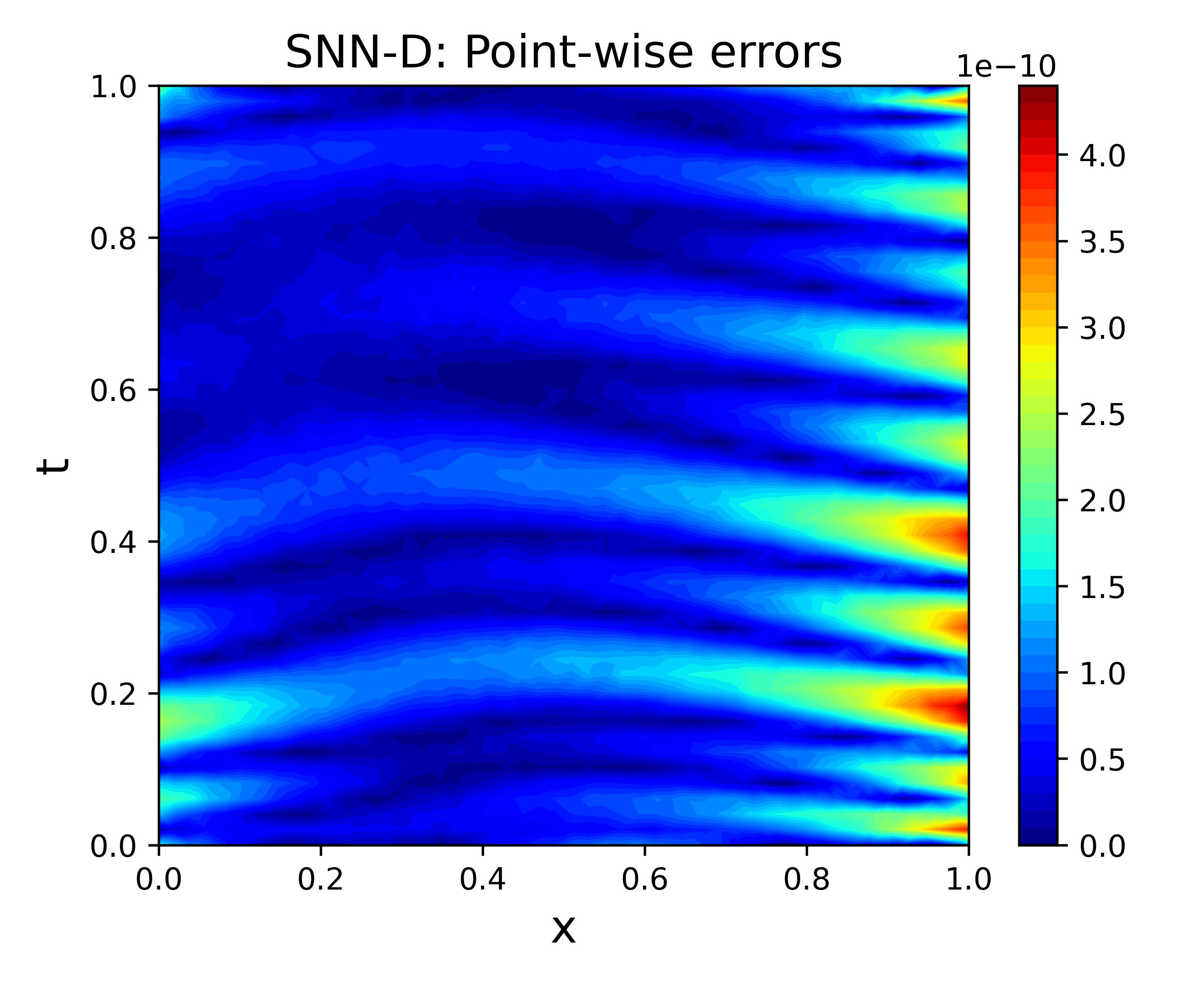}}
\hspace{0.3cm}  
\subcaptionbox{SNN-I}{\includegraphics[width=0.45\linewidth]{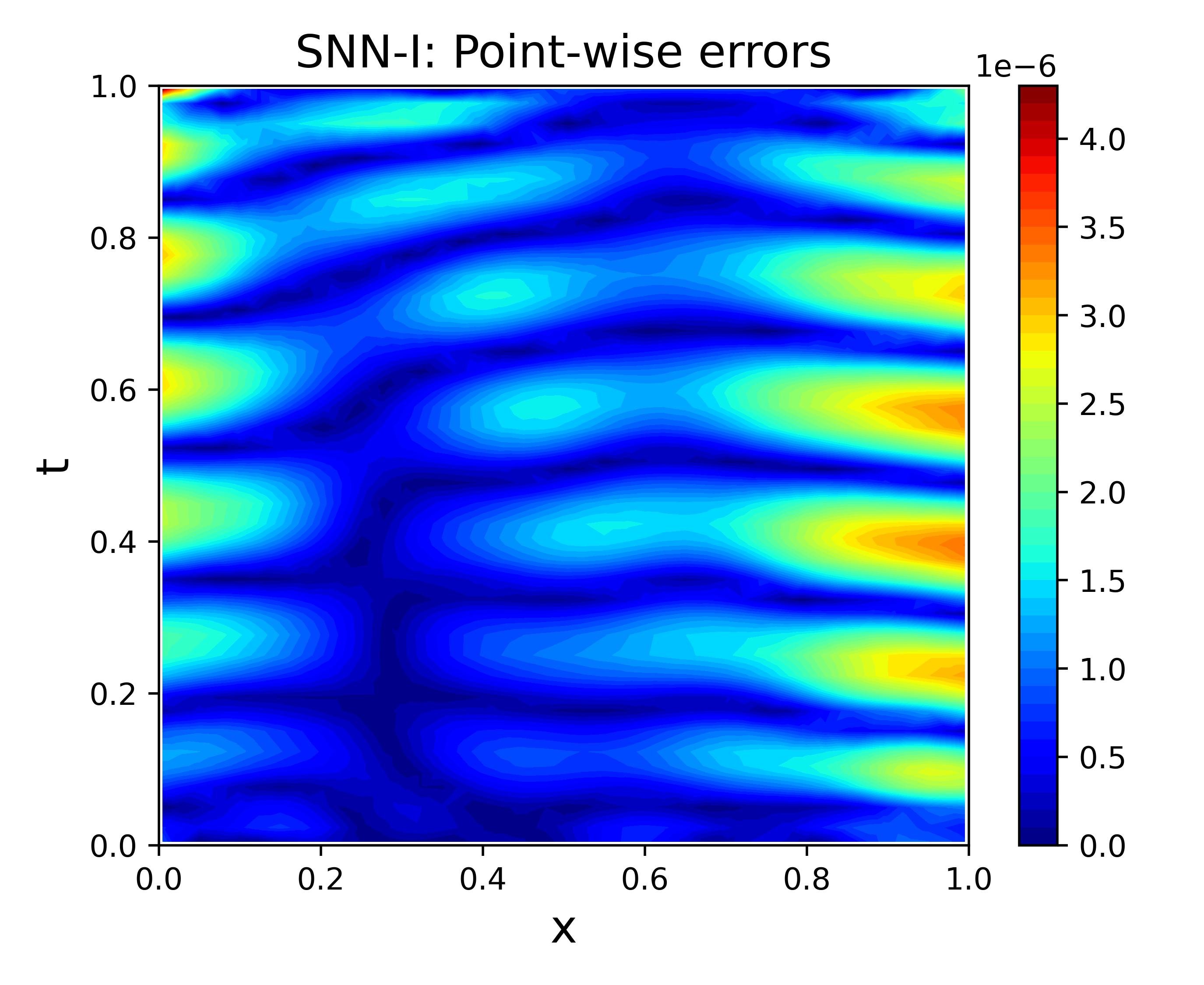}}
\caption{  Point-wise errors of SNN-D and SNN-I for the parabolic equation.}
\label{df_point_error}
\end{figure}

Figure \ref{df_point_error} illustrates the point-wise errors of SNN-D and SNN-I for solving the parabolic equation. Table \ref{df_error} presents the relative \(L_2\) errors of  SNN-D and SNN-I. 
The error of SNN-D is 1.04e-10, and the training epochs of SNN-D is 92.
The error of SNN-I is 1.30e-06, and the training epochs of SNN-I is 136.
This indicates that our methods continue to maintain accuracy and efficiency for solving the parabolic equation.

\begin{table}[!htbp]
\caption{The errors and epochs of SNN-D and SNN-I for the parabolic equation.} 
\label{df_error}
\begin{center}
\small
\begin{tabular}{@{} lccc @{}}
\hline
Method & $\|e\|_{L^2}$ & $\|e\|_{L^\infty}$ & epochs \\
\hline
SNN-D & 1.04e-10 & 4.36e-09 & 92 \\
SNN-I & 1.30e-06  & 4.24e-06 & 136 \\
\hline
\end{tabular}
\end{center}
\end{table}

\subsection{Anisotropic diffusion equation}
\label{Anisotropic diffusion equation}
Consider the following two-dimensional   diffusion equation with anisotropic diffusion coefficient on the domain \(\Omega=(0,1) \times(0,1)\),
\begin{eqnarray}
    \label{Anisotropic_eq}
    -\nabla \cdot(\kappa(x, y) \nabla u)=f(x, y), & (x, y) \in \Omega, \\ 
    u(x, y)=g(x, y), & (x, y) \in \partial \Omega,
\end{eqnarray}
where 
\begin{eqnarray*}
\kappa(x, y)=\left(\begin{array}{cc}
k_1 & 0 \\
0 & k_2
\end{array}\right)
\end{eqnarray*}
We choose  $f (x,y)$ and $g (x,y)$ such that Eq. (\ref{Anisotropic_eq}) has the solution
\(
u(x) = \sin\left(\pi x\right) \sin\left(\pi y \right).
\)

For SNN-D, \(32 \times 32\) points are uniformly sampled throughout the domain, ensuring an uniform distribution of 32 points along each boundary. For SNN-I, a two-dimensional composite Gaussian quadrature formula is used, dividing each dimension into 8 subintervals, with 4 points per subinterval. This strategy leads to 1024 sampling points in the interior and 128 sampling points on the boundary.

Table \ref{ani_error} presents the relative \(L_2\) errors of various methods for solving the anisotropic diffusion equation with different anisotropy ratios. We can see that after 50000 training epochs, PINN and DGM can only maintain certain accuracy when the anisotropy is not very strong, and lose the accuracy when the anisotropic is very strong. However, even when the anisotropy ratio reaches ($1:10^6$), SNN-D achieves the accuracy of 3.92e-09, and SNN-I achieves the accuracy of 2.70e-05. This indicates that SNN-D and SNN-I still can achieve high accuracy for the  diffusion equation with very strong anisotropy ratios.

\begin{table}[!htbp]
\caption{The errors and epochs of PINN, DGM, SNN-D, and SNN-I for the anisotropic diffusion equation with different anisotropy ratios.}\label{ani_error}
\begin{center}
\small
\begin{tabular}{ccccc}
\hline
Method & $k_1:k_2$    & $\|e\|_{L^2}$     & $\|e\|_{L^\infty}$              & epochs \\
\hline
       & $1:10^0$     & 9.54e-03      & 7.33e-03      & 50000      \\
       & $1:10^1$     & 3.14e-02      & 3.14e-02      & 50000      \\
       & $1:10^2$     & 7.41e-01      & 6.47e-01      & 50000      \\
PINN   & $1:10^3$     & 1.68e0        & 1.28e0        & 50000      \\
       & $1:10^4$     & 1.03e01       & 5.6e0         & 50000      \\
       & $1:10^5$     & 9.62e01       & 5.23e0        & 50000      \\
       & $1:10^6$     & 1.04e01       & 5.61e0        & 50000      \\
\hline
       & $1:10^0$     & 5.37e-10      & 1.81e-09      & 276    \\
       & $1:10^1$     & 7.76e-10      & 2.21e-09      & 232    \\
       & $1:10^2$     & 1.64e-09      & 2.88e-09      & 202    \\
SNN-D  & $1:10^3$     & 2.67e-09      & 4.50e-09      & 196    \\
       & $1:10^4$     & 3.90e-09      & 5.30e-09      & 195    \\
       & $1:10^5$     & 3.87e-09      & 4.89e-09      & 195    \\
       & $1:10^6$     & 3.92e-09      & 5.15e-09      & 195    \\
\hline
       & $1:10^0$     & 1.81e-03      & 3.26e-03      & 50000      \\
       & $1:10^1$     & 2.76e-02      & 4.39e-02      & 50000      \\
       & $1:10^2$     & 7.55e-01      & 7.15e-01      & 50000      \\
DGM    & $1:10^3$     & 1.91e0        & 1.45e0        & 50000      \\
       & $1:10^4$     & 1.00e01       & 5.61e0        & 50000      \\
       & $1:10^5$     & 1.02e01       & 5.66e0        & 50000      \\
       & $1:10^6$     & 1.02e01       & 5.66e0        & 50000      \\
\hline
       & $1:10^0$     & 2.97e-06      & 7.97e-06      & 275    \\
       & $1:10^1$     & 4.85e-06      & 1.26e-05      & 211    \\
       & $1:10^2$     & 6.34e-06      & 9.80e-06      & 177    \\
SNN-I  & $1:10^3$     & 2.42e-05      & 2.20e-05      & 170    \\
       & $1:10^4$     & 2.87e-05      & 3.33e-05      & 169    \\
       & $1:10^5$     & 2.72e-05      & 3.27e-05      & 169    \\
       & $1:10^6$     & 2.70e-05      & 3.25e-05      & 169    \\
\hline
\end{tabular}
\end{center}
\end{table}

Table \ref{ani_error2} presents the relative \(L_2\) errors of SNN-D and SNN-I for solving the anisotropic diffusion equation with different strong anisotropy ratios. We can see that our methods can still achieve high accuracy for these tests. 

\begin{table}[!htbp]
\caption{The errors and epochs of SNN-D and SNN-I for the anisotropic diffusion equation with different strong anisotropy ratios.}\label{ani_error2}
\begin{center}
\small
\begin{tabular}{ccccc}\hline
Method & $k_1:k_2$ & $\|e\|_{L^2}$ & $\|e\|_{L^\infty}$ & epochs\\
\hline
      & $10^{0}:10^{-6}$ & 1.46e-09    &   2.38e-09   &    221   \\
      & $10^{1}:10^{-5}$ & 1.99e-09    &   3.01e-09   &    221   \\
      & $10^{2}:10^{-4}$ & 2.00e-09    &   3.27e-09   &    221   \\
SNN-D & $10^{3}:10^{-3}$ & 1.43e-09    &   2.18e-09   &    221   \\
      & $10^{4}:10^{-2}$ & 1.37e-09    &   2.23e-09   &    221   \\
      & $10^{5}:10^{-1}$ & 1.96e-09    &   3.09e-09   &    221   \\
      & $10^{6}:10^{0}$  & 1.96e-09    &   3.11e-09   &    221   \\
\hline
      & $10^{0}:10^{-6}$ & 6.74e-06    &   1.07e-05   &    209   \\
      & $10^{1}:10^{-5}$ & 6.74e-06    &   1.08e-05   &    209   \\
      & $10^{2}:10^{-4}$ & 6.74e-06    &   1.07e-05   &    209   \\
SNN-I & $10^{3}:10^{-3}$ & 6.80e-06    &   1.08e-05   &    209   \\
      & $10^{4}:10^{-2}$ & 6.65e-06    &   1.06e-05   &    209   \\
      & $10^{5}:10^{-1}$ & 6.62e-06    &   1.07e-05   &    209   \\
      & $10^{6}:10^{0} $ & 6.73e-06    &   1.06e-05   &    209   \\
\hline
\end{tabular}
\end{center}
\end{table}

\begin{figure}[htbp]
\centering
\subcaptionbox{SNN-D}{\includegraphics[width=0.45\linewidth]{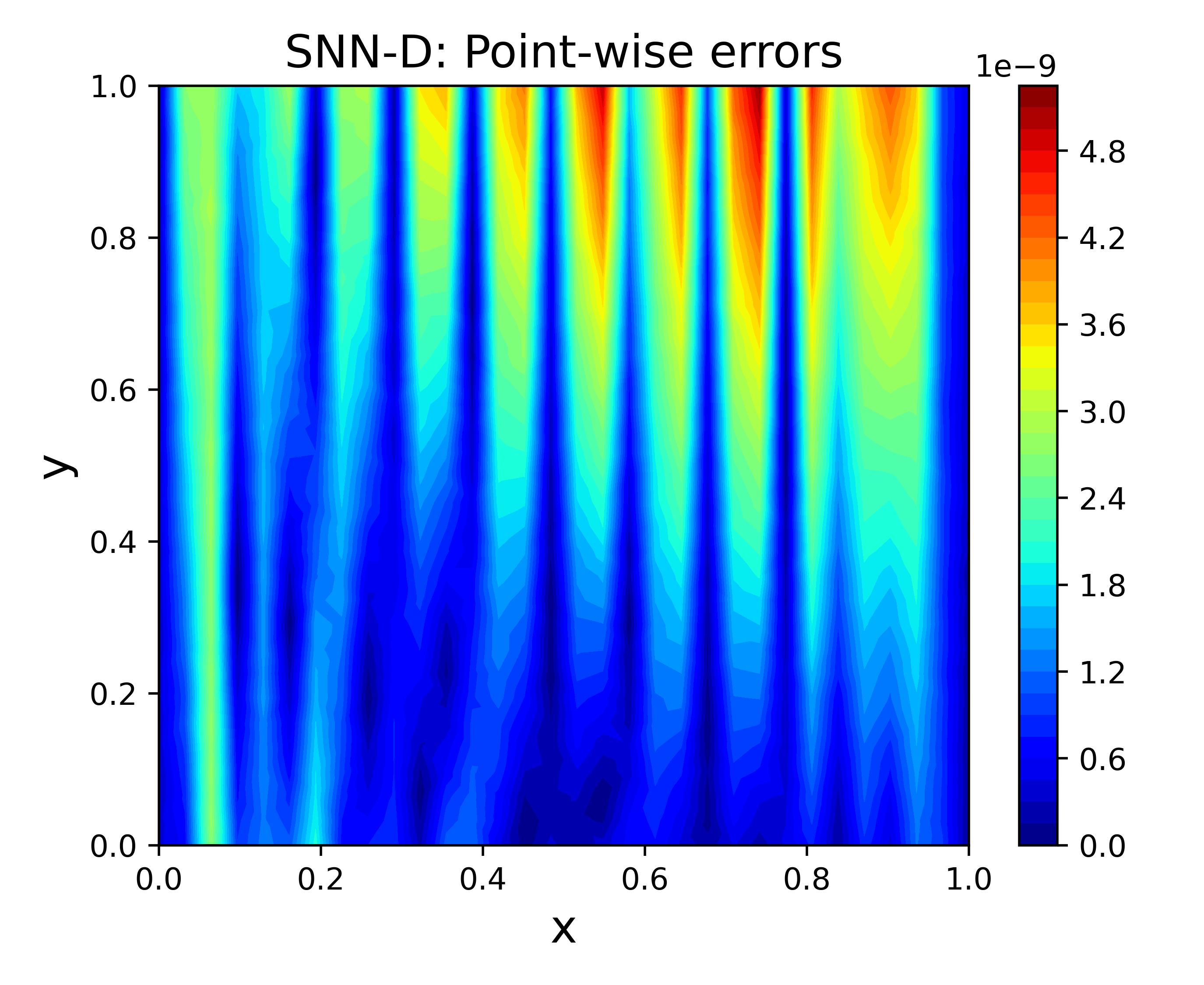}}
\hspace{0.3cm}  
\subcaptionbox{SNN-I}{\includegraphics[width=0.45\linewidth]{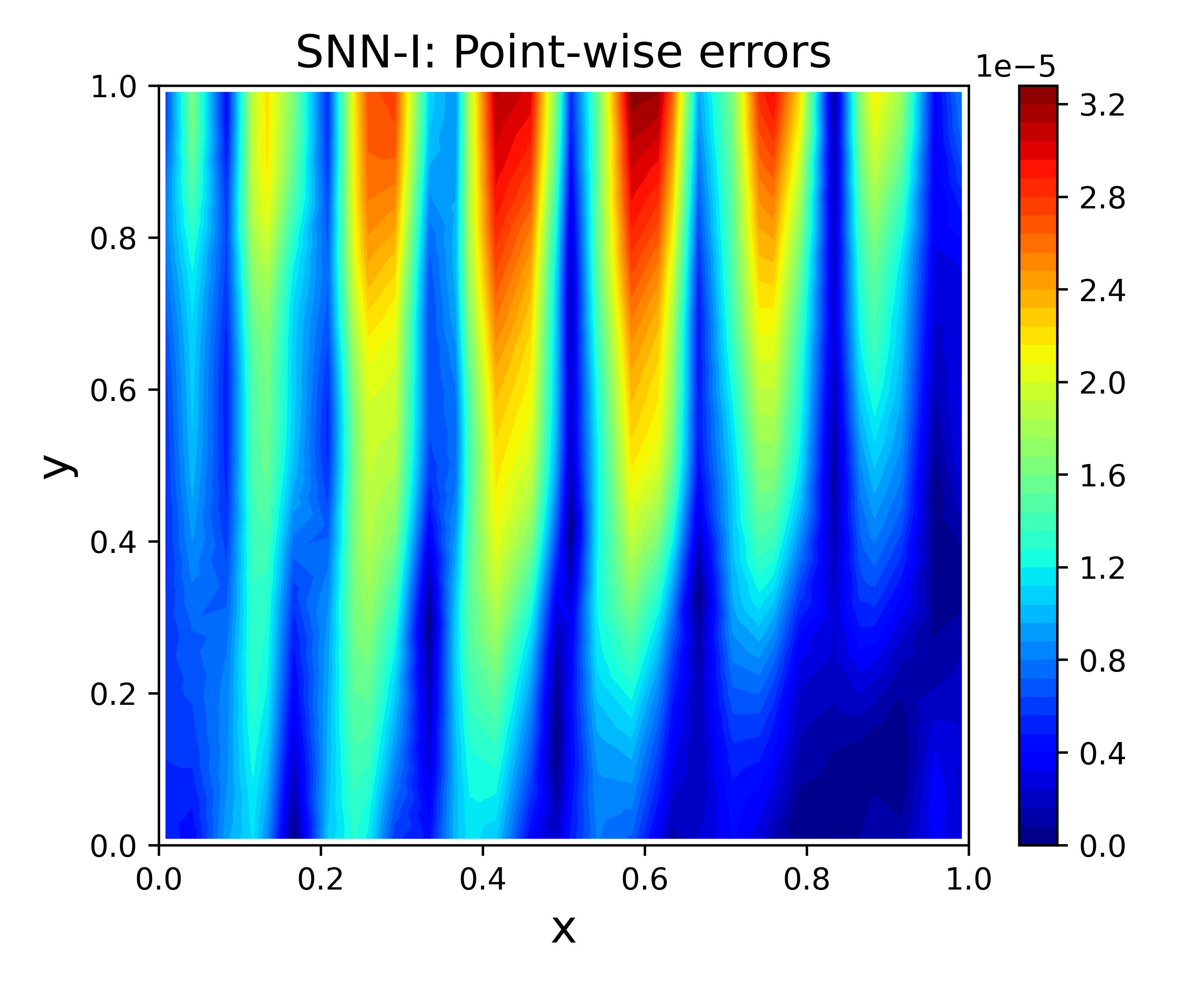}}
\caption{ Point-wise error of SNN-D and SNN-I for the anisotropic diffusion equation with the anisotropy ratio  \(1:10^6\).}
\label{ani_point_error}
\end{figure}

Figure \ref{ani_point_error} illustrates the point-wise errors  of SNN-D and SNN-I for a strong anisotropy ratio (\(1:10^6\)). Hence, SNN-D and SNN-I can maintain accuracy even when the anisotropy is very strong. These tests demonstrate the robustness and adaptability of our method.

\section{Conclusion}
In this paper, we present a subspace method based on neural networks  for solving the partial differential equation with high accuracy. The basic idea of our method is to use some functions based on neural networks  as base functions to span a subspace, then find an approximate solution in this subspace.
We give a general frame of SNN, and design two special algorithms including SNN-D and SNN-I.
Our method is free of parameters  that need to be artificially adjusted. 
The accuracy of SNN is insensitive to the initial parameters of neural networks.  With better initial parameters, SNN requires only a few dozen to several hundred training epochs to rapidly reach the stopping condition. With poorer initial parameters, SNN will adaptively increase the number of training epochs to train these base functions of subspace such that the subspace has effective approximate capability to the solution space of the equation.
Our method is a high-precision deep learning method, the  error can even fall below the level of $10^{-10}$. The cost of training these base functions  of subspace is low, and only one hundred to two thousand epochs are needed for most tests.
The performance of SNN significantly surpasses that of PINN and DGM in terms of the accuracy and computational cost.

\section*{Acknowledgements} 
This work is partially supported by National Natural Science Foundation of China (12071045) and Fund of National Key Laboratory of Computational Physics.


\end{document}